\documentclass[preprint,sort&compress]{elsarticle}
\usepackage{amssymb,amsthm,amsmath,commath,bbm,psfrag,setspace,color,array,booktabs,multirow,url,float}
\usepackage[FIGTOPCAP,nooneline]{subfigure}
\usepackage{algorithm}
\usepackage{algorithmic}

\usepackage{xcolor}

\theoremstyle{definition}
\newtheorem{definition}{Definition}[section]

\newtheorem{theorem}{Theorem}[section]

\newtheorem{lemma}[theorem]{Lemma}

\renewcommand{\baselinestretch}{1.2}
\begin{document}

\begin{frontmatter}

\title{Controlling coexisting attractors of a class of non-autonomous dynamical systems}

\author[EXT]{Zhi Zhang}
\ead{zz326@exeter.ac.uk}

\author[ESPOL,CFD]{Joseph P\'aez Ch\'avez }
\ead{ jpaez@espol.edu.ec}

\author[EXT]{Jan Sieber}
\ead{J.Sieber@exeter.ac.uk}

\author[EXT]{Yang Liu \corref{mycorrespondingauthor}}

\address[EXT]{College of Engineering, Mathematics and
Physical Sciences, University of Exeter, North Park Rd, Exeter EX4
4QF, UK}
\address[ESPOL]{Center for Applied Dynamical Systems and Computational Methods (CADSCOM), Faculty of
Natural Sciences and Mathematics, Escuela Superior Polit\'ecnica del
Litoral, P.O. Box 09-01-5863, Guayaquil, Ecuador}
\address[CFD]{Center for Dynamics, Department of Mathematics, TU Dresden, D-01062 Dresden, Germany}
\cortext[mycorrespondingauthor]{Corresponding author. Tel:
+44(0)1392-724654, e-mail: y.liu2@exeter.ac.uk.}

\begin{abstract}
This paper studies a control method for switching stable coexisting attractors of a class of non-autonomous dynamical systems. The central idea is to introduce a continuous path for the system's trajectory to transition from its original undesired stable attractor to a desired one by varying one of the system parameters according to the information of the desired attractor. The behaviour of the control is studied numerically for both non-smooth and smooth dynamical systems, using a soft-impact and a Duffing oscillators as examples. Special attention is given to identify the regions where the proposed control strategy is applicable by using path-following methods implemented via the continuation platform COCO. It is shown that the proposed control concept can be implemented through either using an external control input or varying a system parameter. Finally, extensive numerical results are presented to validate the proposed control methods.
\end{abstract}


\end{frontmatter}

\section{Introduction}
Non-autonomous dynamical systems exhibit a rich variety of different
long-term behaviours coexisting for a given set of parameters, which
is referred as multistability or coexisting attractors. It is a
common phenomenon in science and nature that relies crucially on the
initial conditions of the system. For example, multistability can be
found in many engineering systems, such as the vibro-impact capsule
\cite{liu2017capsule,liu2020bifurcation}, the electronic circuits
\cite{arecchi1982hopping, beasley1983comment}, the gas laser
\cite{arecchi1982experimental} and the drilling system
\cite{pavlovskaia2015modelling,liu2017drilling}. It is also a
fundamental property of biological systems, such as the spiking
neurons \cite{foss1996multistability}, tumor progression
\cite{huang2013genetic}, cell fate transitions
\cite{suel2006excitable,huang2005cell} and cell cycle control
\cite{yao2011origin,battogtokh2004bifurcation}. Control of
multistability has been an active research field in the past few
decades, see \cite{pisarchik2014control}. These coexisting
attractors are extremely sensitive to noise due to their fractally
interwoven basins of attraction. Switching the coexisting attractors
while protecting against perturbation-induced basin hopping is vital
in practice, particularly for engineering applications. Since the
engineering systems may present different performance at different
coexisting attractors (e.g. energy efficiency), control the
switching among these stable attractors could offer them more
flexibility. For example, stick-slip motion and constant rotation
coexist in the drilling system
\cite{liu2019drilling,liu2020drilling}, and prevention of the
stick-slip motion can maintain the drilling at a high penetration
rate so reducing the operation cost. For this reason, numerous
methods have been developed to control multistability, and a control
strategy for avoiding perturbation-induced attractor switching is
always desirable.

In the present work, we will develop a novel control method for
controlling the coexisting attractors of a class of non-autonomous
dynamical systems, including both smooth and non-smooth dynamical
systems.
It is well known that the OGY method
\cite{edward1990controlling} was initially developed for stabilizing
an unstable periodic orbit embedded in a chaotic motion via
adjusting the system parameter in a small neighbourhood. For the
same purpose, Pyragas \cite{pyragas1992continuous} developed a
delayed feedback controller to achieve the stabilisation of unstable
periodic orbits in a chaotic system. In
\cite{pyragas1993experimental}, Pyragas and Tama{\v{s}}evi{\v{c}}ius
used the delayed feedback control method to stabilise an analogue
circuit. In \cite{pyragas2018act,pyragas2019state}, Pyragas and
Pyragas introduced the act-and-wait concept to reduce the dimension
of phase space of the systems with the delayed feedback control.
Although the above studies focused on stabilising unstable periodic
attractors embedded in the chaotic attractor, those control concepts
also can be applied to control coexisting attractors. For example,
Lai \cite{lai1996driving} constructed a hierarchy of paths for
targeting the desired attractor by introducing a feedback
perturbation to the system with fractal basins. Wang \emph{et al}.
\cite{wang2016geometrical} adopted the similar concept to control an
undesired attractor to a desired one for a low dimensional network
system. In \cite{zhang2020calculating}, Zhang \emph{et al}. used the
delayed feedback control to switch the unwanted attractors to a
period-1 attractor for a soft impacting oscillator, and in
\cite{chavez2020numerical}, P{\'a}ez Ch{\'a}vez \emph{et al}.
developed a numerical approach for analysing the dynamical
properties of the non-smooth systems with time delays. On the other
hand, there are also many other methods for controlling the
multistability of dynamical systems. Arecchi \emph{et al}.
\cite{arecchi1985generalized} discovered that the external noise can
bridge the coexisting states of a forced Duffing equation. In
\cite{pisarchik2000annihilation}, Pisarchik and Goswami employed a
slow external periodic perturbation to the system parameter of a
bistable system to annihilate its coexisting attractors. In
\cite{liu2013intermittent}, an intermittent control was designed to
provide an impulsive force for non-autonomous dynamical systems to
switch between coexisting attractors, which was verified by both
numerical and experimental results. Liu and P{\'a}ez Ch{\'a}vez
\cite{liu2017controlling} developed a linear augmentation control
law to control the multistability of a soft impacting oscillator,
and analysed the dynamical properties of the control law by using
the path-following (continuation) techniques for non-smooth
dynamical systems. Nevertheless, it should be noted that most of
these control methods for attractor's switching of non-autonomous
dynamical systems relied on external input, while few works have
concerned about the switching by utilising system's properties, e.g.
modulation of system parameter. It is obvious that the latter
approach is less invasive, so the original system could be
maintained, which is advantageous to some systems whose external
input is hard to access, such as the machining process
\cite{stepan2017,stepan2019,yan2021basins}. Therefore, in the
present work, we will focus on the control method that relies on the
original properties of the system only to achieve the switching
among coexisting attractors.

According to \cite{pisarchik2014control}, Pisarchik and Feudel
classified the existing methods for controlling multistability into
three categories: feedback, non-feedback and stochastic controls.
They have suggested that the most efficient way of ensuring a
predefined behaviour for the system is to annihilate all the other
coexisting attractors. However, annihilation of undesired attractors
in some dynamical system could change the existing structure of
solutions leading to the emergence of new complex basins of
attraction. One of the simplest ways to achieve this could be to
apply an impulsive external perturbation to direct system's
trajectory from one basin to another one
\cite{pisarchik2014control}. In \cite{kaneko1990}, Kaneko's study
suggested that a short pulse can cause system state to jump from one
attractor to another one if the pulse's amplitude is sufficiently
large. Chizhevsky \emph{et al}. proposed to apply a short-pulsed
perturbation in the form that switching the system off and on for a
very short time, so the system can run from a different initial
condition. The problem of these methods is that the short pulse was
applied in the form of non-feedback. If the amplitude of the pulse
is small, the system regains the same attractor after a few
iterations. In addition, these methods are only effective for the
systems with few coexisting attractors, and the short-pulsed control
may become uncertain if more attractors coexist with fractal basins.
Hence, continuous control by utilising system's feedback signals is
more preferable for non-smooth dynamical systems that have
considerable numbers of coexisting attractors, e.g. in the
near-grazing dynamics \cite{yin2019complex}.

In the present work, the primary focus is to address the
continuous switching between two of coexisting stable
attractors by varying a system parameter without affecting their
original dynamics. In order to achieve this, we propose a continuous control method that can adjust a system parameter based on the information of trajectory of the desired
attractor. This control method is applied to the controlled
system continuously until its trajectory is sufficiently close to the desired one. We develop two control strategies based on this concept, the so-called linear and nonlinear control strategies, where the former is implemented through an external control input and the latter is applied via a system parameter. The advantage of the nonlinear control strategy is that it depends only on the original properties of system parameter and does not rely on any external input. To demonstrate the applicability of the proposed control method to both non-smooth and smooth systems, an impact and a Duffing oscillators were employed in the present paper.

The rest of this paper is organized as follows. Section
\ref{sec-mathpre} introduces the mathematical model of the
periodically forced mechanical oscillator subjected to a one-sided
soft constraint and some mathematical preliminaries for the proposed control method. In
Section \ref{sec-numerical}, the effectiveness of this control
method on a single-degree-of-freedom oscillator with piecewise-smooth nonlinearity is studied. In Section
\ref{sec-smooth}, the proposed control method is adopted to the Duffing
oscillator. Finally, concluding remarks are drawn in Section \ref{sec-conclusion}.

\section{Design of feedback control strategy}
  \label{sec-mathpre}
We consider periodically forced systems with a single control input of the form
\begin{equation}\label{main-prob-u}
\begin{split}
\dot{Y}_\mathrm{u}(\tau)&=F(\tau,Y_\mathrm{u}(\tau),u(\tau)), \\
Y_\mathrm{u}(\tau_{0})&=Y_{\mathrm{u},0},
\end{split}
\end{equation}
where the input $u(\tau)$ is a scalar function of time, and $\dot{Y}$ denotes
differentiation with respect to time $\tau$. We assume that
the uncontrolled system, \eqref{main-prob-u} with $u(\tau)=0$, has an
attractor $Y_\mathrm{d}(\cdot)$ (where the subscript in $Y_\mathrm{d}$
stands for ``desired''), but that the initial condition
$Y_{\mathrm{u},0}$ is away from $Y_\mathrm{d}(0)$, and possibly
outside the basin of attraction of $Y_\mathrm{d}(\cdot)$.
We also assume that this ``desired'' attractor $Y_d(\cdot)$ is (internally) stable in the sense that it has no positive Lyapunov exponents. In practice our method is intended to be applied to periodic orbits.
Our test examples will be single-degree-of-freedom oscillators with
periodic forcing and multiple coexisting attractors. The
typical scenario we envisage is that $Y_{\mathrm{u},0}$ is on one of the other (``undesirable'') attractors of the uncontrolled system.

\subsection{Single-degree-of-freedom oscillator with piecewise-smooth nonlinearity}
\label{sec:soft-impact}
\begin{figure}[h!]
\centering
\includegraphics[width=3.1in]{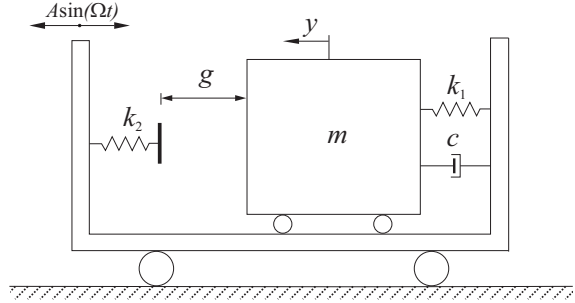}
\caption{Physical model of the single-degree-of-freedom oscillator with piecewise-smooth nonlinearity.}\label{fig-model}
\end{figure}

In this section, we will use the single-degree-of-freedom oscillator with piecewise-smooth nonlinearity shown in Fig.~\ref{fig-model} as an example to study the proposed control
method. Soft impacts occur in mechanical systems when an object hits
an obstacle with a negligible mass but with a non-negligible
stiffness, see e.g. \cite{bernardo2008piecewise, de2008control,
lazarek2020detection, serdukova2021post, makarenkov2012dynamics}. As
can be seen from Fig.~\ref{fig-model}, it is assumed that the
discontinuity boundary is fixed at $x=e$, with $e>0$ being the
nondimensional gap. The equations of motion of the oscillator
are in form \eqref{main-prob-u}, where $Y(\tau):=(x(\tau),v(\tau))^{T}$. We will consider three different cases. Defining for scalar $\tilde{u}$
\begin{align}\label{main-prob}
F_\mathrm{gen}(Y,\tilde{u},u_e):=\left[\begin{array}{cc}
  0 & I  \\
-I -\beta H(x-e) & -2\zeta \\
\end{array}\right]Y+\left[\begin{array}{cc}
  0 \\
\beta (e+u_e)H(x-e-u_e)+\tilde{u}
\end{array}\right],
\end{align}
the three cases for control input are
\begin{align}\label{flin}
  F_\mathrm{lin}(\tau,Y,u)&:=F_\mathrm{gen}(Y,a\omega^2\sin(\omega\tau) + u,0)\mbox{,}\\
  \label{famp}
  F_a(\tau,Y,u)&:=F_\mathrm{frc}(Y,(a+u)\omega^2\sin(\omega\tau),0)\mbox{,}\\
  \label{fgap}
  F_e(\tau,Y,u)&:=F_\mathrm{frc}(Y,a\omega^2\sin(\omega\tau),u)\mbox{.}
\end{align}
The first case, $F_\mathrm{lin}$ has linear control input and the second case, $F_a$, has parametric control input, varying the forcing amplitude, while the third case adjusts the gap $e$.
The function $H(\cdot)$ stands for the Heaviside step function. In
the right-hand side $F_\mathrm{frc}$, defined in \eqref{main-prob}, the variables (including time) and
parameters of the system are nondimensionalised according to
\begin{equation*}\label{eq-nondim}
\setstretch{2.25}\begin{array}{r@{}lcr@{}lcr@{}lcr@{}l}
\omega_{n}=&\mbox{ }\sqrt{\dfrac{k_1}{m}}, & & \tau=&\mbox{ }\omega_{n}t, & &
\omega=&\mbox{ }\dfrac{\Omega}{\omega_{n}}, & & \zeta=&\mbox{ }\dfrac{c}{2m\omega_{n}},\\
x=&\mbox{ }\dfrac{y}{y_{0}}, & & e=&\mbox{ }\dfrac{g}{y_{0}}, & &
a=&\mbox{ }\dfrac{A}{y_{0}}, & & \beta=&\mbox{
}\dfrac{k_{2}}{k_{1}},
\end{array}
\end{equation*}
where $y_{0}>0$ is an arbitrary reference length, $\omega_n$ is the
natural frequency, $\omega$ is the frequency ratio, $\beta$ is the
stiffness ratio, $\zeta$ is the damping ratio, $e$ is the
nondimensional gap between the mass and the secondary spring and $a$
is the nondimensionalised amplitude of the external excitation.

The right-hand side given in Eq.~\eqref{main-prob} is a typical
non-autonomous dynamical system with the one-sided elastic constraint
considered as the non-smoothness that can lead to complex phenomena,
such as the grazing bifurcation \cite{yin2019complex,yin20}, the
coexistence of multistable attractors
\cite{ing2007experimental,ing2010bifurcation} and the chaotic motions
\cite{lazarek2020detection}. Here, we consider the situation that
system \eqref{main-prob-u} with $u(\tau)=0$ has many coexisting
attractors within some specific ranges of the system parameters.
In particular, we consider a stable attractor $Y_\mathrm{d}$ with non-positive Lyapunov exponents exists, which means this attractor is periodic or quasi-periodic.
Usually, this $Y_\mathrm{d}(\cdot)$, our desired attractor, is a stable periodic orbit. For $Y_\mathrm{d}$, small changes in system parameters do not affect the response of the system
significantly, such that the attractor $Y_\mathrm{d}$ will persist.

\subsection{Distance-reducing feedback control}
\label{sec:control}
Let us introduce the following definition.
\begin{definition}\label{def-control}
For a dynamical system $\dot{Y}(\tau)=F(\tau,Y(\tau),0)$, $Y_\mathrm{u}(\tau_0)=Y_{\mathrm{u},0}$ (so, of type \eqref{main-prob-u} with $u=0$) with two stable coexisting attractors $Y_\mathrm{c}(\cdot)$ and $Y_\mathrm{d}(\cdot)$, if
there exists a continuous control $u(\tau)$ such that system \eqref{main-prob-u} with this control $u(\cdot)$ and $Y_{\mathrm{u},0}=Y_\mathrm{c}(\tau_0)$ satisfies $Y_\mathrm{u}(\tau)-Y_\mathrm{d}(\tau)\to0$ for $\tau\to\infty$, we say that $Y_\mathrm{c}$ is controllable to $Y_\mathrm{d}$ (by $u$). If $Y_\mathrm{u}(\tau^*)=Y_\mathrm{d}(\tau^*)$ for some finite time $\tau^*$ we say that $Y_\mathrm{c}$ is controllable to $Y_\mathrm{d}$ in finite time.
\end{definition}

Since $Y_\mathrm{d}(\cdot)$ is an attractor for $u(\tau)=0$, the system will follow $Y_\mathrm{d}(\tau)$ for $\tau>\tau^*$, if we set $u(\tau)=0$ for $\tau>\tau^*$ after controlling to $Y_d$ in finite time $\tau^*$.
Let us assume that system \eqref{main-prob-u} is controllable, and
there exists a control input $u(\cdot)$ that can be used for
controlling its multistability.

Define the difference between the desired and the current states as
\begin{align*}
  d(\tau)&:=Y_\mathrm{d}(\tau)-Y_\mathrm{u}(\tau)\mbox{,\quad such that}\\
  \Delta(\tau)&:=\langle d(\tau), d(\tau) \rangle\geq0
\end{align*}
is the distance at time $\tau$ from the desired attractor
$Y_\mathrm{d}$. We will investigate a simple feedback control strategy
that aims to reduce this distance $\Delta$ over time, such that
$[\mathrm{d}/\mathrm{d}\tau]\Delta(\tau)<0$. Requiring that $u(\tau)$ has permissible values between $-M_1$ and $M_1$ and Lipschitz constant $M_2$, we adjust $u$ at every time $\tau\geq\tau_0$ according to
\begin{align}
  \label{u:design}
  \dot{u}(\tau)&=M_2\operatorname{sign}\left[-\frac{\partial}{\partial u}\frac{\mathrm{d}}{\mathrm{d}\tau}\Delta(\tau)\right]H(M_1+u)H(M_1-u)\mbox{,}
\end{align}
starting from $u(\tau_0)=0$, for
as long as $\Delta(\tau)>\epsilon$ with some tolerance $\epsilon\ll1$. In Eq.~\eqref{u:design}, $H(M_1-u)H(M_1+u)$ is the indicator function for the interval $[-M_1,M_1]$ (recall that $H$ is the Heaviside function). The term $[\mathrm{d}/\mathrm{d}\tau]\Delta(\tau)$ is a function of $(\tau,Y_\mathrm{u}(\tau),Y_\mathrm{d}(\tau),u(\tau))$, making the derivative with respect to $u$ non-zero:
\begin{align}
  \label{dist:form}
  \frac{\partial}{\partial u}\left[\frac{\mathrm{d}}{\mathrm{d}\tau}\Delta\right](\tau,Y_\mathrm{u},Y_\mathrm{d},u)&=
  2\langle Y_\mathrm{d}-Y_\mathrm{u},\partial_uF(\tau,Y_\mathrm{u},u)\rangle\mbox{.}
\end{align}
In the practical algorithms described in Sections \ref{external_subsection} and \ref{sec:parcont}, we will apply Eq.~\eqref{u:design} only when $[\mathrm{d}/\mathrm{d}\tau]\Delta(\tau)>0$, otherwise, we will set $\dot{u}(\tau)=0$.

\subsection{Implementation with finite sampling step --- linear case}\label{external_subsection}
Let us assume the sampling time step, $h>0$, such that
$\tau_i=\tau_0+ih$. If the control input $u$ enters the right-hand
side linearly with a constant coefficient vector $b$, such that
$F(\tau,Y,u)=F(\tau,Y)+bu$ (as in example $F_\mathrm{lin}$ given in
\eqref{main-prob}, where $b=(0,1)^\mathrm{T}$), the notation in the
definition \eqref{u:design} of the control $u$ and the resulting
expression \eqref{dist:form} simplifies. Thus, we formulate Algorithm~1 for this common case separately. In our formulation of Algorithm~1 with the sampling step $h$, the control input $u(\tau)$ is of type ``zero-order hold''. That is, $u(\tau)$ is a constant $u_i$ on the sampling interval $[\tau_i,\tau_{i+1}]$, and the Lipschitz constant $M_2$ applies to changes per time step: $|u_{i+1}-u_i|\leq hM_2$.

\renewcommand{\baselinestretch}{2.0}
\begin{algorithm}\label{linear_algo}
\caption{Linear control}
\begin{algorithmic}
\STATE \textbf {Step 0}:

Choose $M_{1}$ and $M_{2}$, where $M_{1}$ is the boundary of
$|u(\tau)|$, $M_{2}$ is the boundary of $|\dot{u}(\tau)|$, and
$\dot{u}(\tau):=\tfrac{\mathrm{d}u(\tau)}{\mathrm{d}\tau}$. Take the
initial control $u(\tau_{0})=0\in M_{1}$, and set the iteration index
$i=0$ and the time step $h$. Below we denote $U(\tau)=bu(\tau)$. \WHILE{the termination criterion
$\langle d(\tau_\mathrm{i}), d(\tau_\mathrm{i}) \rangle\leq
\epsilon$ is not satisfied,}

\STATE \textbf {Step 1}:
Compute the range of $\dot{u}(\tau_\mathrm{i})$ that satisfies
the following criterion
\begin{align}\label{condition}
\langle d(\tau_\mathrm{i}), \dot{U}(\tau_\mathrm{i}) \rangle h\ge
&2\langle d(\tau_\mathrm{i}), F(\tau_\mathrm{i}, Y_{\mathrm{d}}(\tau_\mathrm{i}))-
F(\tau_\mathrm{i}, Y_\mathrm{u}(\tau_\mathrm{i})) -U(\tau_\mathrm{i})\rangle \\
&+\langle F(\tau_\mathrm{i},
Y_{\mathrm{d}}(\tau_\mathrm{i}))-F(\tau_\mathrm{i},
Y_\mathrm{u}(\tau_\mathrm{i})) -U(\tau_\mathrm{i}),  F(\tau_\mathrm{i},
Y_{\mathrm{d}}(\tau_\mathrm{i}))-F(\tau_\mathrm{i},
Y_\mathrm{u}(\tau_\mathrm{i})) -U(\tau_\mathrm{i}) \rangle h \nonumber,
\end{align}
calling this range the feasible range for $\dot u(\tau_i)$.
\STATE \textbf {Step 2}:
If this feasible range is greater than $M_{2}$ and ensure that $u_{i+1}$
satisfies $M_{1}$, take the minimum value of
$|\dot{u}(\tau_\mathrm{i})|$ and go to \textbf {Step 3}. Otherwise,
choose the value $\dot{u}(\tau_\mathrm{i})$ such that both $M_{1}$
and $M_{2}$ are satisfied and is the closest to the feasible range
of $\dot{u}(\tau_\mathrm{i})$, and go to \textbf{Step 3}.
(NB. Use $\dot{u}(\tau_\mathrm{i})=0$ if $d(\tau_\mathrm{i})$ decreases at $\tau_\mathrm{i}$.)

\textbf {Step 3}:

Use $u_{i+1}=u_{i}+\dot{u}(\tau_\mathrm{i})h$
for input $u$ in Eq.~\eqref{main-prob-u} for $\tau>\tau_i$.
Increase $i$ by 1 and return to \textbf {Step 1}.

\ENDWHILE
\end{algorithmic}
\end{algorithm}

\renewcommand{\baselinestretch}{1.2}

After introducing the algorithm, the following main theorem can be
obtained.
\begin{theorem}\label{Main_1}
Let the dynamical system
  $\dot
  Y(\tau)=F(\tau,Y)+(\mathbf{0},u(\tau))^\mathrm{T}$
  with $u=0$ have two stable coexisting attractors $Y_\mathrm{c}$ and
  $Y_\mathrm{d}$, and let $M_{1,2}$ be bounded intervals in
  $\mathbb{R}$. We assume that there exists a time
  $\tau^{\ast}>\tau_{0}$, such that for all sufficiently small
  sampling steps $h$ the control $U$ and trajectory $Y_\mathrm{u}$ obtained by
  Algorithm~1 with bounds $M_{1,2}$ satisfy
\begin{align}\label{inequality_control}
&\big|\langle d(\tau_{0}), d(\tau_{0}) \rangle
+\sum_{i=0}^{n^{\ast}}\big[2\langle d(\tau_{i}),
F(\tau_{i}, Y_{\mathrm{d}}(\tau_{i}))-F(\tau_{i},
Y_\mathrm{u}(\tau_{i})) -U(\tau_{i})\rangle h-
\langle d(\tau_{i}), \dot{U}(\tau_{i}) \rangle h^2\nonumber\\
&+\langle F(\tau_{i},
Y_{\mathrm{d}}(\tau_{i}))-F(\tau_{i}, Y_\mathrm{u}(\tau_{i}))
-U(\tau_{i}),  F(\tau_{i},
Y_{\mathrm{d}}(\tau_{i}))-F(\tau_{i}, Y_\mathrm{u}(\tau_{i}))
-U(\tau_{i}) \rangle h^2\big]\big|\leq c_{1} h^2
\end{align}
(here $d(\tau)=Y_\mathrm{d}(\tau)-Y_\mathrm{u}(\tau)$ and the $\tau_i$ are from
the interval partition
$[\tau_{0},\tau^{\ast})=\bigcup_{i=0}^{n^{\ast}}[\tau_\mathrm{i},\tau_{i+1})$
with $\tau_{i+1}=\tau_\mathrm{i}+h$ for $i=0,\ldots,n^{\ast}$).  Then
these two stable attractors are controllable by this external
controller.
\end{theorem}

According to \textbf{Theorem \ref{Main_1}}, the following lemma can
be obtained.
\begin{lemma}
For any two of stable coexisting attractors of system
(\ref{main-prob}), if the control sequence  $\{u_{i}\}$,
$i=1,\cdots,n^{*}$ generated by Algorithm 1 satisfies
$u(\tau_{0})=0$, $u(\tau)\in M_{1}$, $\dot{u}(\tau)\in M_{2}$, where
$M_{1,2}$ are the bounded intervals of $\mathbb{R}$, inequality
(\ref{inequality_control}) and $| u(\tau^{*}) |\leq c_{1} h$, then
these two stable attractors are controllable by this external
control sequence.
\end{lemma}

\subsection{Implementation with finite sampling step --- nonlinear case}
\label{sec:parcont}
By adjusting the accessible parameter of non-autonomous dynamical
systems, such as the amplitude of excitation, system energy is
altered, so the dynamical property of the system can be controlled
(see e.g.
\cite{edward1990controlling,lai1996driving,wang2016geometrical}). In
this subsection, we will apply the control method in Algorithm 1 to
an accessible system parameter to achieve the switching between
two stable coexisting attractors. The detailed algorithm of the new control method (Algorithm 2) is
given as below, and the following theorem is introduced.

\renewcommand{\baselinestretch}{2.0}
\begin{algorithm}\label{adaptive_algo}
\caption{Nonlinear control}
\begin{algorithmic}
\STATE \textbf {Step 0}:
Choose $M_\mathrm{p,1}$ and $M_\mathrm{p,2}$, which are the
boundaries of $u_\mathrm{p}(\tau)$ and $\dot{u}_\mathrm{p}(\tau)$,
respectively. Take the initial control $u_\mathrm{p}(\tau_{0})=0\in
M_\mathrm{p,1}$, and set the iteration $i:=0$ and the time step $h$.
\WHILE{the termination criterion $\langle
d_\mathrm{p}(\tau_\mathrm{i}), d_\mathrm{p}(\tau_\mathrm{i})
\rangle\leq \epsilon$ is not satisfied,} \STATE \textbf {Step 1}:

Compute the feasible range of $\dot{u}_\mathrm{p}(\tau_\mathrm{i})$
to satisfy the following criterion
\begin{align}\label{condition_2}
&\langle d_\mathrm{p}(\tau_\mathrm{i}), \frac{\mathrm{D} F(\tau_\mathrm{i},Y_\mathrm{u}
(\tau_\mathrm{i}),u_\mathrm{p}(\tau_\mathrm{i}))}{\mathrm{D} u_\mathrm{p}} \rangle
\dot{u}_\mathrm{p}(\tau_\mathrm{i}) h\nonumber\\
&\qquad\ge 2\langle d_\mathrm{p}(\tau_\mathrm{i}), F(\tau_\mathrm{i}, Y_{\mathrm{d}}
(\tau_\mathrm{i}))-F(\tau_\mathrm{i},Y_\mathrm{u}(\tau_\mathrm{i}),u_\mathrm{p}
(\tau_\mathrm{i}))\rangle -\langle d_\mathrm{p}(\tau_\mathrm{i}), \frac{\mathrm{D}
 F(\tau_\mathrm{i},Y_\mathrm{u}(\tau_\mathrm{i}),u_\mathrm{p}(\tau_\mathrm{i}))}{\mathrm{D}  \tau} \rangle h\nonumber\\
&\qquad\quad +\langle F(\tau_\mathrm{i}, Y_{\mathrm{d}}(\tau_\mathrm{i}))-
F(\tau_\mathrm{i}, Y_\mathrm{u}(\tau_\mathrm{i}),u_\mathrm{p}(\tau_\mathrm{i})),
F(\tau_\mathrm{i}, Y_{\mathrm{d}}(\tau_\mathrm{i}))-F(\tau_\mathrm{i},
Y_\mathrm{u}(\tau_\mathrm{i}),u_\mathrm{p}(\tau_\mathrm{i}))  \rangle h \nonumber\\
&\qquad\quad -\langle d_\mathrm{p}(\tau_\mathrm{i}) ,
\frac{\mathrm{D}
F(\tau_\mathrm{i},Y_\mathrm{u}(\tau_\mathrm{i}),u_\mathrm{p}(\tau_\mathrm{i}))}{\mathrm{D}
x} \rangle h,
\end{align}

\STATE \textbf {Step 2}:

If this range is greater than $M_{2}$ and ensure that
$u_\mathrm{p,i}$ satisfied $M_\mathrm{p,1}$, take the minimum value
of $|\dot{u}_\mathrm{p}(\tau_\mathrm{i})|$ and go to \textbf{Step
3}. Otherwise, choose the value
$\dot{u}_\mathrm{p}(\tau_\mathrm{i})$ such that both
$M_\mathrm{p,1}$ and $M_\mathrm{p,2}$ are satisfied and is closest
to the feasible range of $\dot{u}_\mathrm{p}(\tau_\mathrm{i})$. Then
go to  \textbf{Step 3}.

\textbf {Step 3}:

Use
$u_\mathrm{p,i+1}=u_\mathrm{p,i}+\dot{u}_\mathrm{p}(\tau_\mathrm{i})h$
for Eq.~\eqref{main-prob-u} for $\tau>\tau_i$. Increment $i$ by 1 and return to \textbf{Step 1}.

\ENDWHILE
\end{algorithmic}
\end{algorithm}

\renewcommand{\baselinestretch}{1.2}

\begin{theorem}\label{Main_2}
Let the dynamical system
  $\dot
  Y(\tau)=F(\tau,Y, u_{p}(\tau))$
  with $u_{p}=0$ have two stable coexisting attractors $Y_\mathrm{c}$ and
  $Y_\mathrm{d}$, and let $M_{p,\kappa}$, $\kappa=1, 2$ be bounded intervals in
  $\mathbb{R}$. We assume that there exists a time
  $\tau^{\ast}>\tau_{0}$, such that for all sufficiently small
  sampling steps $h$ the control $u_{p}$ and trajectory $Y_\mathrm{u}$ obtained by
  Algorithm~2 with bounds $M_{p,\kappa}$ satisfy
\begin{align}\label{inequality_adapcontrol}
&\big|\langle d_\mathrm{p}(\tau_{0}), d_\mathrm{p}(\tau_{0}) \rangle
+\sum_{i=0}^{n^{*}}\big[2\langle d_\mathrm{p}(\tau_\mathrm{i}),
F(\tau_\mathrm{i},
Y_{\mathrm{d}}(\tau_\mathrm{i}))-F(\tau_\mathrm{i},
Y_{\mathrm{u}}(\tau_\mathrm{i}),u_\mathrm{p}(\tau_\mathrm{i})) \rangle h-
\langle d_\mathrm{p}(\tau_\mathrm{i}), \frac{\mathrm{D} F(\tau_\mathrm{i},
Y_{\mathrm{u}}(\tau_\mathrm{i}),u_\mathrm{p}(\tau_\mathrm{i}))}{\mathrm{D}  \tau} \rangle h^2\nonumber\\
&+\langle F(\tau_\mathrm{i}, Y_{\mathrm{d}}(\tau_\mathrm{i}))-F(\tau_\mathrm{i},
Y_{\mathrm{u}}(\tau_\mathrm{i}),u_\mathrm{p}(\tau_\mathrm{i})) ,  F(\tau_\mathrm{i},
Y_{\mathrm{d}}(\tau_\mathrm{i}))-F(\tau_\mathrm{i}, Y_{\mathrm{u}}(\tau_\mathrm{i}),
u_\mathrm{p}(\tau_\mathrm{i}))  \rangle h^2\nonumber\\
&-\langle d_\mathrm{p}(\tau_\mathrm{i}) , \frac{\mathrm{D}
F(\tau_\mathrm{i},Y_{\mathrm{u}}(\tau_\mathrm{i}),u_\mathrm{p}(\tau_\mathrm{i}))}{\mathrm{D}
Y} \rangle h^{2}- \langle d_\mathrm{p}(\tau_\mathrm{i}),
\frac{\mathrm{D}
F(\tau_\mathrm{i},Y_{\mathrm{u}}(\tau_\mathrm{i}),u_\mathrm{p}(\tau_\mathrm{i}))}{\mathrm{D}
u_\mathrm{p}} \dot{u}_\mathrm{p}(\tau_\mathrm{i}) \rangle  h^{2}
\big] \big|\leq c_{2} h^2,
\end{align}
(here $d(\tau)=Y_\mathrm{d}(\tau)-Y_\mathrm{u}(\tau)$ and the $\tau_i$ are from
the interval partition
$[\tau_{0},\tau^{\ast})=\bigcup_{i=0}^{n^{\ast}}[\tau_\mathrm{i},\tau_{i+1})$
with $\tau_{i+1}=\tau_\mathrm{i}+h$ for $i=0,\ldots,n^{\ast}$).  Then
these two stable attractors are controllable by this external
controller.
\end{theorem}

According to \textbf{Theorem \ref{Main_2}}, the following lemma can
be obtained.
\begin{lemma}
For any two of stable coexisting attractors of system
(\ref{main-prob}), if the control sequence  $\{u_\mathrm{p,i}\}$,
$i=1,\cdots,n^{*}$ generated by Algorithm 2 satisfies
$u_\mathrm{p}(\tau_{0})=0$, $u_\mathrm{p}(\tau)\in M_\mathrm{p,1}$,
$\dot{u}_\mathrm{p}(\tau)\in M_\mathrm{p,2}$, where $M_\mathrm{p,1}$
and $M_\mathrm{p,2}$ are the bounded intervals of $\mathbb{R}$,
inequality condition \eqref{inequality_adapcontrol} and $|
u_\mathrm{p}(\tau^{*}) |\leq c_{2} h$, then these two stable
attractors are controllable by varying the system parameter $p$.
\end{lemma}

\section{Control of non-smooth dynamical systems}\label{sec-numerical}
\begin{figure}[h!]
\centering
\includegraphics[width=0.55\textwidth]{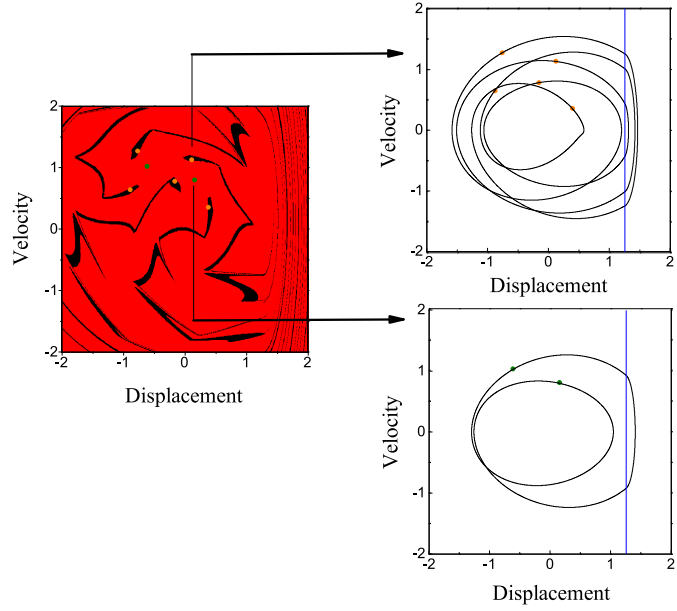}
\caption{(Colour online) Basins of attraction of the impacting
system computed for $\zeta= 0.01, ~e = 1.26, ~a = 0.7, ~\beta = 28 ~
\textrm{and}~ \omega= 0.85$. Orange dots denote the period-$5$
attractor with black basin, and green dots represent the period-$2$
attractor with red basin. The right panels present the trajectories
of the period-$5$ and the period-$2$ attractors on the phase plane,
where blue lines indicate the impact boundary.}\label{coexist}
\end{figure}

This section will show the effectiveness of the proposed control
methods by using the impact oscillator shown in
Fig.~\ref{fig-model}, which is a typical non-smooth dynamical system
exhibiting many coexisting attractors at its near-grazing dynamics
\cite{ing2010bifurcation}. The following parameters were used for
which two stable attractors, a period-$2$ and a period-$5$
responses, coexist as shown in Fig.~\ref{coexist}.
\begin{equation*}\label{parameter}
\zeta= 0.01, ~e = 1.26, ~a = 0.7, ~\beta = 28 ~ \textrm{and}~
\omega= 0.85.
\end{equation*}

\subsection{Linear control}

\begin{figure}[h!]
\centering
\includegraphics[width=0.85\textwidth]{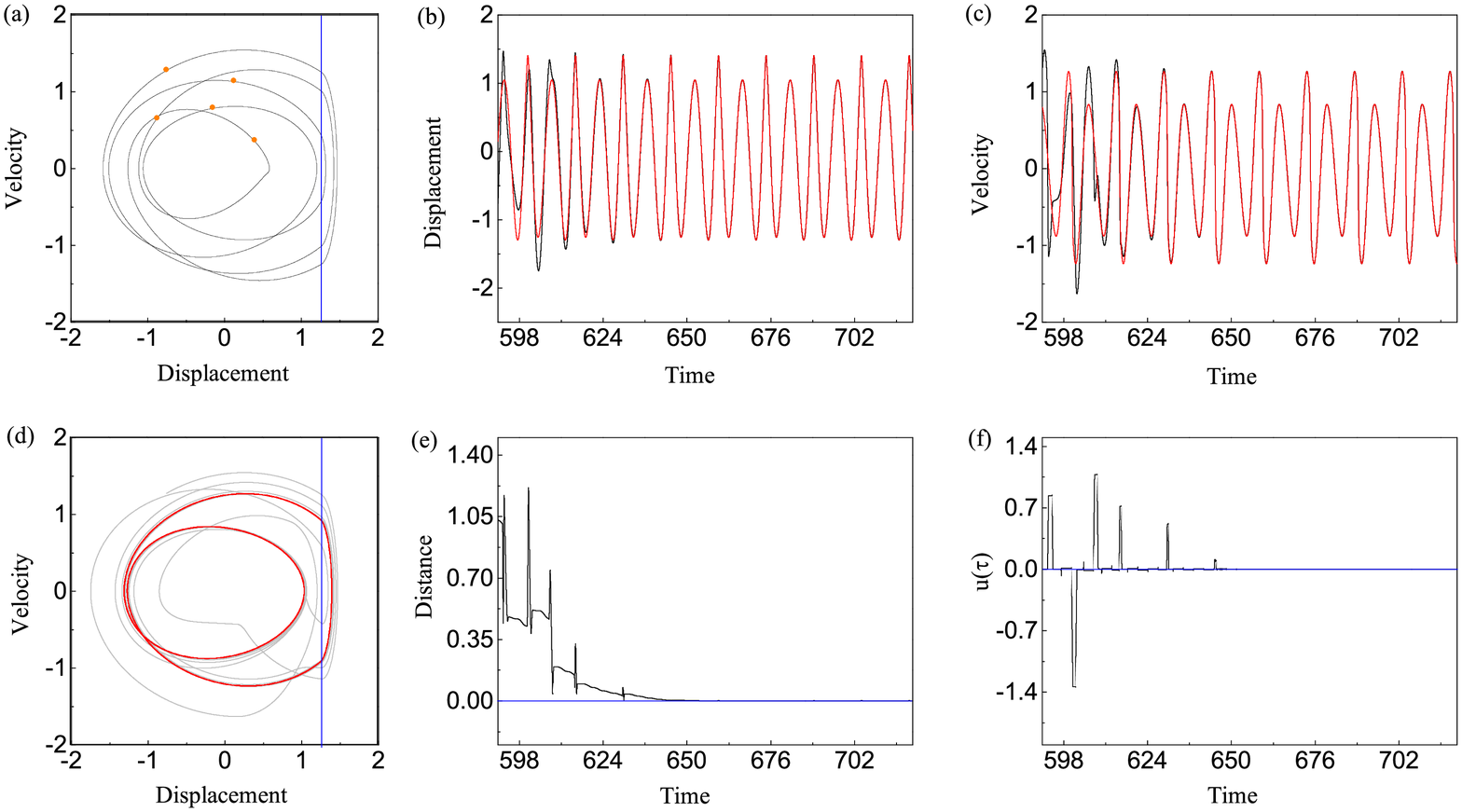}
\caption{(Colour online) (a) The period-5 response on the phase
plane with the Poincar\'{e} sections denoted by orange dots. (b)
Time histories of the desired (red line) and the current (black
line) displacements of the system under the external control
strategy (Algorithm 1) with $M_{1}=5$ and $M_{2}=3$. (c) Time
histories of the desired (red line) and the current (black line)
velocities of the system. (d) Trajectory of the system on the phase
plane under the external control strategy, where grey and red lines
represent the transient and the steady-state responses,
respectively. (e) Time history of the distance between the desired
and the controlled trajectories in $2$-norm. (f) Time history of the
control sequence generated by the external control strategy. Blue
lines in (a) and (d) indicate the impact boundary, while the blue
lines in (e) and (f) mark the zero reference. The result was
computed for $\zeta= 0.01, ~e = 1.26, ~a = 0.7, ~\beta = 28 ~
\textrm{and}~ \omega= 0.85$.} \label{firstcontrol}
\end{figure}

\begin{figure}[h!]
\centering
\includegraphics[width=0.85\textwidth]{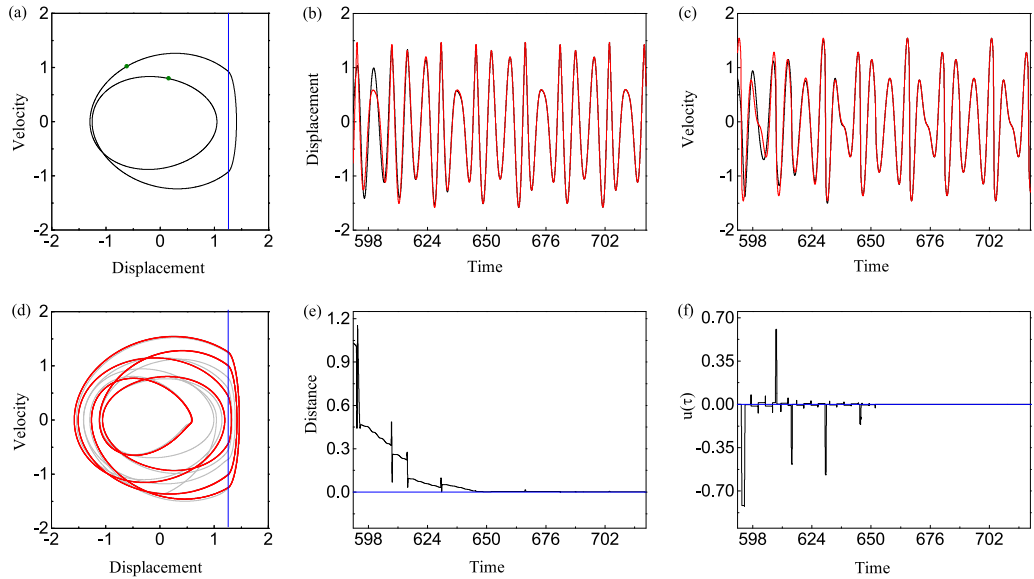}
\caption{(Colour online) (a) The period-2 response on the phase
plane with the Poincar\'{e} sections denoted by green dots. (b) Time
histories of the desired (red line) and the current (black line)
displacements of the system under the external control strategy
(Algorithm 1) with $M_{1}=5$ and $M_{2}=3$. (c) Time histories of
the desired (red line) and the current (black line) velocities of
the system. (d) Trajectory of the system on the phase plane under
the external control strategy, where grey and red lines represent
the transient and the steady-state responses, respectively. (e) Time
history of the distance between the desired and the controlled
trajectories in $2$-norm. (f) Time history of the control sequence
generated by the external control strategy. Blue lines in (a) and
(d) indicate the impact boundary, while the blue lines in (e) and
(f) mark the zero reference. The result was computed for $\zeta=
0.01, ~e = 1.26, ~a = 0.7, ~\beta = 28 ~ \textrm{and}~ \omega=
0.85$.} \label{firstdifferent}
\end{figure}

Firstly, the dynamical response of the impacting system of type
(\ref{main-prob}) with right-hand side $F_\mathrm{lin}$ given in \eqref{flin} is presented in Fig.~\ref{firstcontrol} at where
the original period-5 attractor was switched to the period-2
attractor by using the external control strategy. The time step of
the simulations was fixed at $h=0.002$, and the control strategy was
implemented at $\tau=591.358$. As can be seen from
Figs.~\ref{firstcontrol}(b) and (c), the period-5 response
experienced a transition and was settled down to the period-2
response around $\tau=637$. Fig.~\ref{firstcontrol}(d) shows the
transition by grey line on the phase plane and indicates the
steady-state response by red line. Figs.~\ref{firstcontrol}(e) and
(f) also demonstrate the effectiveness of the control where the
distance between the present and the desired attractors was reduced
once the control sequence was applied. It can be seen that the
overall trend of the distance was decreased, and according to the
simulation, it was about $0.026$ at $\tau=637$ and was about nil
after $\tau=650$. In addition, it can be seen from
Fig.~\ref{firstcontrol}(f) that, during the control process, no
control was applied when the distance between the two trajectories
was decreasing. Then a continuous increase in the control signal to
$0.832$ at $\tau=593.01$ and a continuous decrease to $0$ at
$\tau=594.65$ were recorded. Thereafter, the control experienced
intermittent control actions, and the amplitudes of the control
actions decreased as the two trajectories were closer. Finally, the
control was turned off when the control target was achieved.

Fig.~\ref{firstdifferent} shows the control result from the period-2
to the period-5 attractor by using the external control strategy.
According to the simulation, the control was switched on at
$\tau=591.358$, and the distance between the two trajectories was
decreased to $0.01$ at $\tau=643.5$. The control signal became zero
gradually when the distance was sufficiently small at $\tau=651.9$.

\subsection{Nonlinear control}
Next, we will test the effectiveness of the nonlinear control strategy (Algorithm 2) by varying the amplitude of excitation $a$ and the gap $e$, as given in \eqref{famp} and \eqref{fgap}.
Again, the control target here is to
switch the response of system \eqref{main-prob-u} with right-hand side \eqref{famp} between the
period-5 and period-2 attractors by varying its amplitude of
excitation.

\begin{figure}[h!]
\centering
\includegraphics[width=0.85\textwidth]{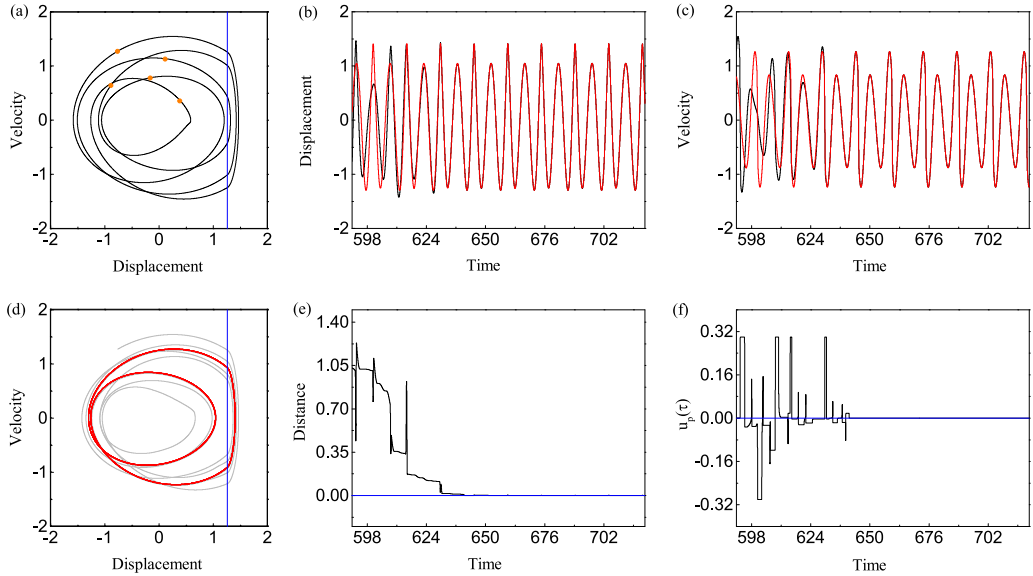}
\caption{(Colour online) (a) The period-5 response on the phase
plane with the Poincar\'{e} sections denoted by orange dots. (b)
Time histories of the desired (red line) and the current (black
line) displacements of the system under the nonlinear
control strategy (Algorithm 2) by varying the amplitude of
excitation with $M_\mathrm{a,1}=0.3$ and $M_\mathrm{a,2}=5$. (c)
Time histories of the desired (red line) and the current (black
line) velocities of the system. (d) Trajectory of the system on the
phase plane under the nonlinear control strategy, where grey
and red lines represent the transient and the steady-state
responses, respectively. (e) Time history of the distance between
the desired and the controlled trajectories in $2$-norm. (f) Time
history of the control sequence generated by the nonlinear
control strategy. Blue lines in (a) and (d) indicate the impact
boundary, while the blue lines in (e) and (f) mark the zero
reference. The result was computed for $\zeta= 0.01, ~e = 1.26, ~a =
0.7, ~\beta = 28 ~ \textrm{and}~ \omega= 0.85$.} \label{second3}
\end{figure}

Fig.~\ref{second3} shows the control of the impacting system from
the period-5 attractor to the period-2 attractor by varying its
amplitude of excitation $a$. Based on our calculation, when the
nonlinear control strategy was applied at $\tau=591.358$,
the controlled trajectory of the system experienced a transition as
observed in Figs.~\ref{second3}(b), (c) and (d). As can be seen from
Fig.~\ref{second3}(e), the distance between the desired and the
controlled trajectories in 2-norm was decreased indicating the
controlled trajectory approached to the desired one, and this
distance was reduced to $0.001$ at $\tau \approx 640.64$.
Fig.~\ref{second3}(f) presents the time history of the control
signal $u_\mathrm{p}(\tau)$ that initiated from $\tau=591.358$ and
terminated at $\tau=640.854$. Thereafter, $u_\mathrm{p}(\tau)=0$,
and the amplitude of excitation was back to its original value,
$a=0.7$.

\begin{figure}[h!]
\centering
\includegraphics[width=0.85\textwidth]{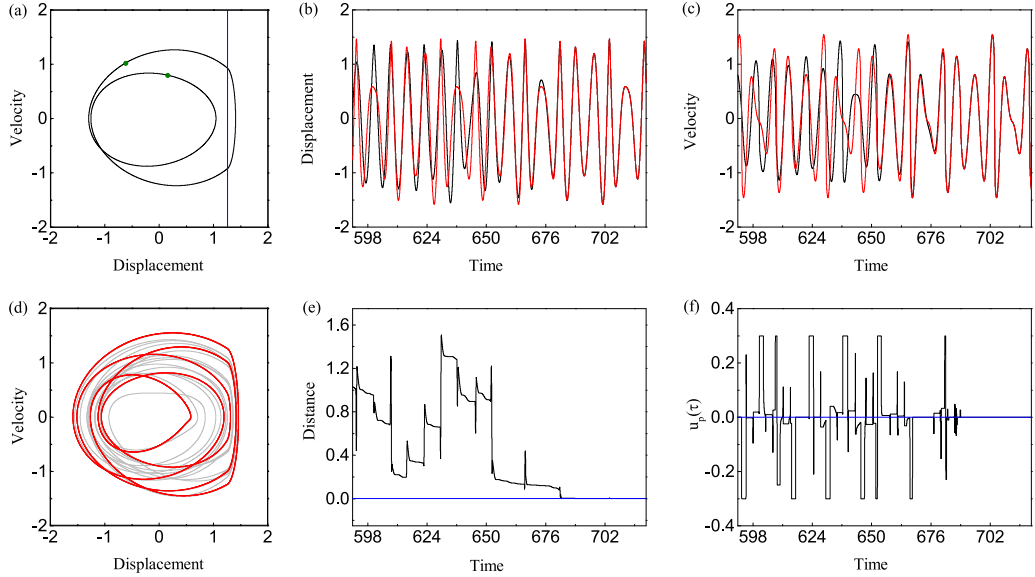}
\caption{(Colour online) (a) The period-2 response on the phase
plane with the Poincar\'{e} sections denoted by green dots. (b) Time
histories of the desired (red line) and the current (black line)
displacements of the system under the nonlinear control
strategy (Algorithm 2) by varying the amplitude of excitation with
$M_\mathrm{a,1}=0.3$ and $M_\mathrm{a,2}=5$. (c) Time histories of
the desired (red line) and the current (black line) velocities of
the system. (d) Trajectory of the system on the phase plane under
the nonlinear control strategy, where grey and red lines
represent the transient and the steady-state responses,
respectively. (e) Time history of the distance between the desired
and the controlled trajectories in $2$-norm. (f) Time history of the
control sequence generated by the nonlinear control
strategy. Blue lines in (a) and (d) indicate the impact boundary,
while the blue lines in (e) and (f) mark the zero reference. The
result was computed for $\zeta= 0.01, ~e = 1.26, ~a = 0.7, ~\beta =
28 ~ \textrm{and}~ \omega= 0.85$.} \label{seconddifferent3}
\end{figure}

The control from the period-2 attractor to the period-5 attractor by
varying the amplitude of excitation $a$ is presented in
Fig.~\ref{seconddifferent3}. The control strategy was applied at
$\tau=591.358$, but the distance between the two trajectories was
not decreased continuously. Therefore, compared to the external
control strategy, it took a longer time for the system to settle
down to the period-5 attractor. According to the simulation, the
control was switched off at $\tau=689.1$ when the distance was
reduced to $0.001$.

When applying the control input $u_\mathrm{p}(\tau)$ to the gap
$e$, as given in \eqref{fgap}, the control aims again to switch the
response of the soft-impact system \eqref{main-prob-u} with right-hand side $F_e(\tau,Y,u_\mathrm{p})$ given in \eqref{fgap}, between the period-5 and period-2 attractors by varying its gap $e$.

\begin{figure}[h!]
\centering
\includegraphics[width=0.82\textwidth]{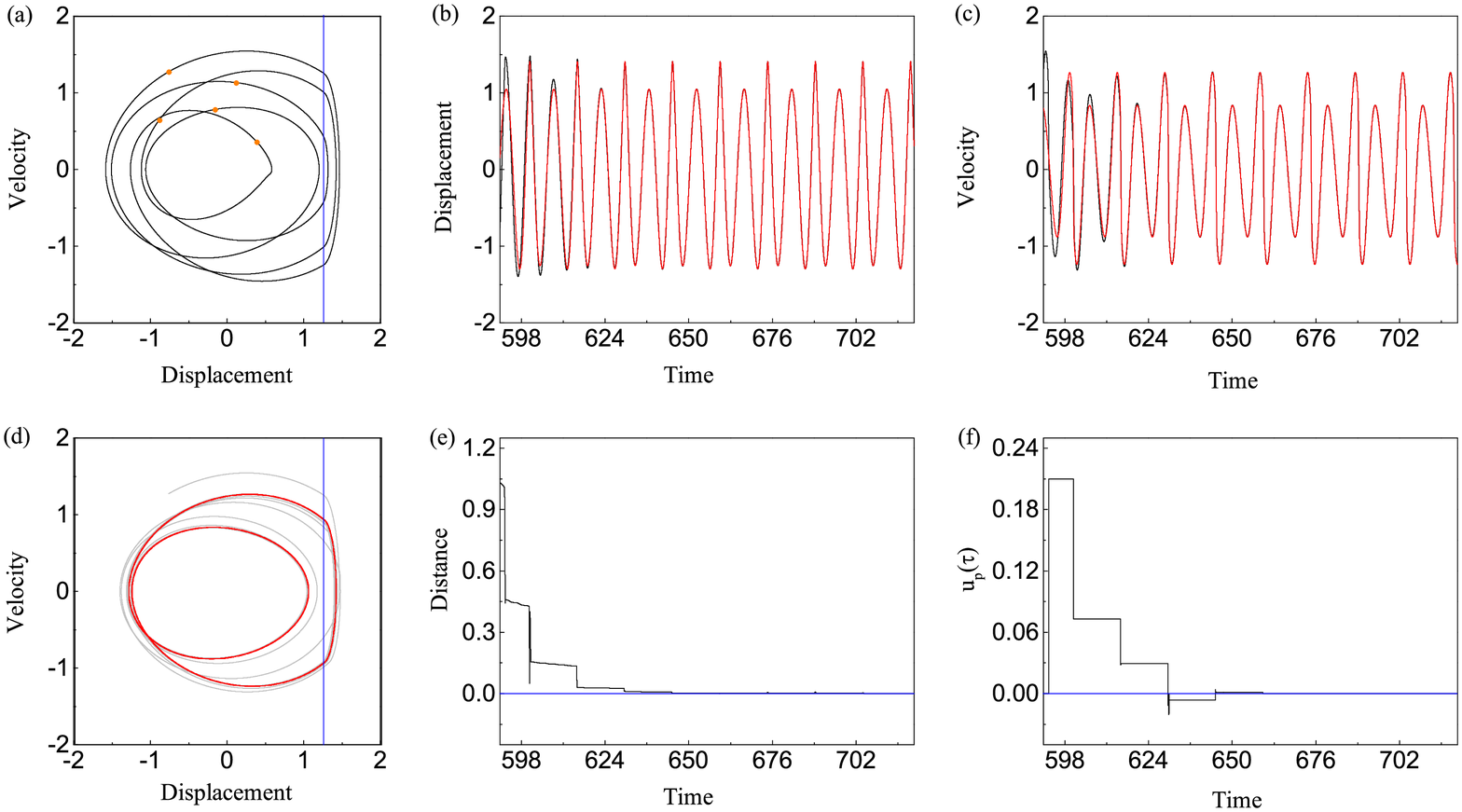}
\caption{(Colour online) (a) The period-5 response on the phase
plane with the Poincar\'{e} sections denoted by orange dots. (b)
Time histories of the desired (red line) and the current (black
line) displacements of the system under the nonlinear
control strategy (Algorithm 2) by varying system's gap with
$M_\mathrm{e,1}=0.3$ and $M_\mathrm{e,2}=5$. (c) Time histories of
the desired (red line) and the current (black line) velocities of
the system. (d) Trajectory of the system on the phase plane under
the nonlinear control strategy, where grey and red lines
represent the transient and the steady-state responses,
respectively. (e) Time history of the distance between the desired
and the controlled trajectories in $2$-norm. (f) Time history of the
control sequence generated by the nonlinear control
strategy. Blue lines in (a) and (d) indicate the impact boundary,
while the blue lines in (e) and (f) mark the zero reference. The
result was computed for $\zeta= 0.01, ~e = 1.26, ~a = 0.7, ~\beta =
28 ~ \textrm{and}~ \omega= 0.85$.} \label{third1}
\end{figure}

\begin{figure}[h!]
\centering
\includegraphics[width=0.82\textwidth]{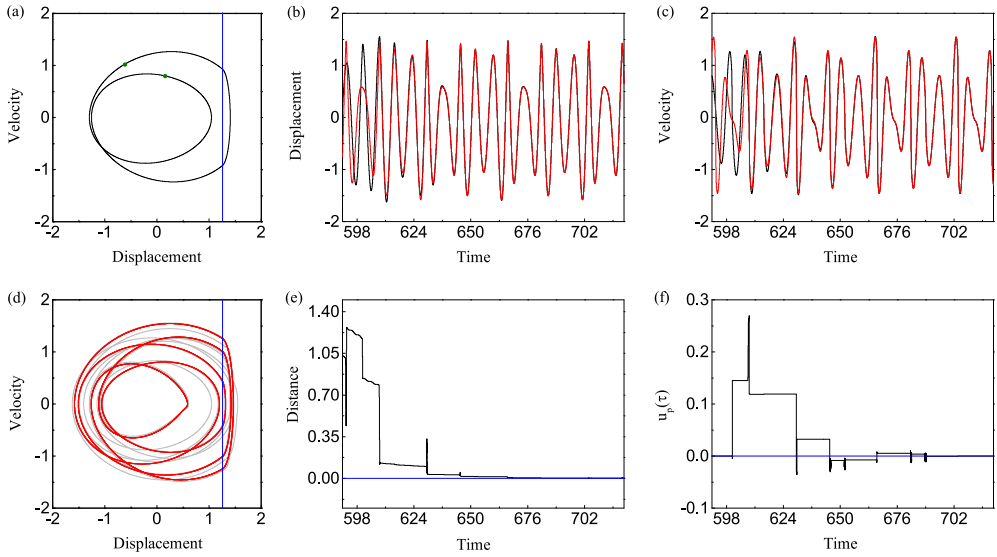}
\caption{(Colour online) (a) The period-2 response on the phase
plane with the Poincar\'{e} sections denoted by green dots. (b) Time
histories of the desired (red line) and the current (black line)
displacements of the system under the nonlinear control
strategy (Algorithm 2) by varying system's gap with
$M_\mathrm{e,1}=0.3$ and $M_\mathrm{e,2}=5$. (c) Time histories of
the desired (red line) and the current (black line) velocities of
the system. (d) Trajectory of the system on the phase plane under
the nonlinear control strategy, where grey and red lines
represent the transient and the steady-state responses,
respectively. (e) Time history of the distance between the desired
and the controlled trajectories in $2$-norm. (f) Time history of the
control sequence generated by the nonlinear control
strategy. Blue lines in (a) and (d) indicate the impact boundary,
while the blue lines in (e) and (f) mark the zero reference. The
result was computed for $\zeta= 0.01, ~e = 1.26, ~a = 0.7, ~\beta =
28 ~ \textrm{and}~ \omega= 0.85$.} \label{thirddifferent1}
\end{figure}

Fig.~\ref{third1} presents the control result from the period-5 to
the period-2 attractor by varying system's gap $e$. The control
strategy was applied at $\tau=591.358$ and was switched off at
$\tau=659.6$. During the control period, the trajectory of the
system experienced a transition, and the distance between the
controlled and the desired trajectories was decreased from $1.03$ to
$0.001$. Within the same time duration, the control signal
$u_\mathrm{p}(\tau)$ reached a maximum value
$u_\mathrm{p}(\tau)=0.21$ and decreased to $-0.02$ at
$\tau=630.326$. Thereafter, the control signal did not change
significantly and reduced to nil gradually after $\tau=659.49$.

To demonstrate the switching from the period-2 to the period-5
attractor by varying system's gap, Fig.~\ref{thirddifferent1}
presents the control result. Based on the calculation, the control
strategy was applied at $\tau=591.358$, and the distance between the
two trajectories was reduced from $1.03$ to $0.001$ at
$\tau=703.794$. Thereafter, the control signal $u_\mathrm{p}(\tau)$
decreased to nil, and the control target was achieved. Compared to
the control result shown in Fig.~\ref{third1}, the transition from
the period-2 to the period-5 attractor took a longer time and had a
more complex transient response, which was due to the complexity of
the period-5 response.

\subsection{Bifurcation analysis of the coexisting
attractors}\label{sec-coex-IO}

In this section we will study in detail the effect of the control
parameters (excitation amplitude $a$ and mass-spring gap $e$) on the
period-2 and period-5 coexisting attractors studied in the previous
section. For this purpose we will employ path-following methods for
piecewise-smooth dynamical systems, using the continuation platform
COCO \cite{dankowicz2013recipes}. The precise COCO-implementation
for the impact oscillator \eqref{main-prob} can be found in a
previous publication by the authors \cite{liu2017controlling}, which
will be adopted in the present work.

\begin{figure}[h!]
\centering \psfrag{T}{\normalsize Contact
time}\psfrag{X}{\normalsize$x(\tau)$}\psfrag{V}{\normalsize$v(\tau)$}\psfrag{a}{\large$a$}
\psfrag{as}{(a)}\psfrag{b}{(b)}\psfrag{c}{(c)}
\psfrag{d}{(d)}\psfrag{e}{(e)}\psfrag{f}{(f)} \psfrag{g}{(g)}
\includegraphics[width=\textwidth]{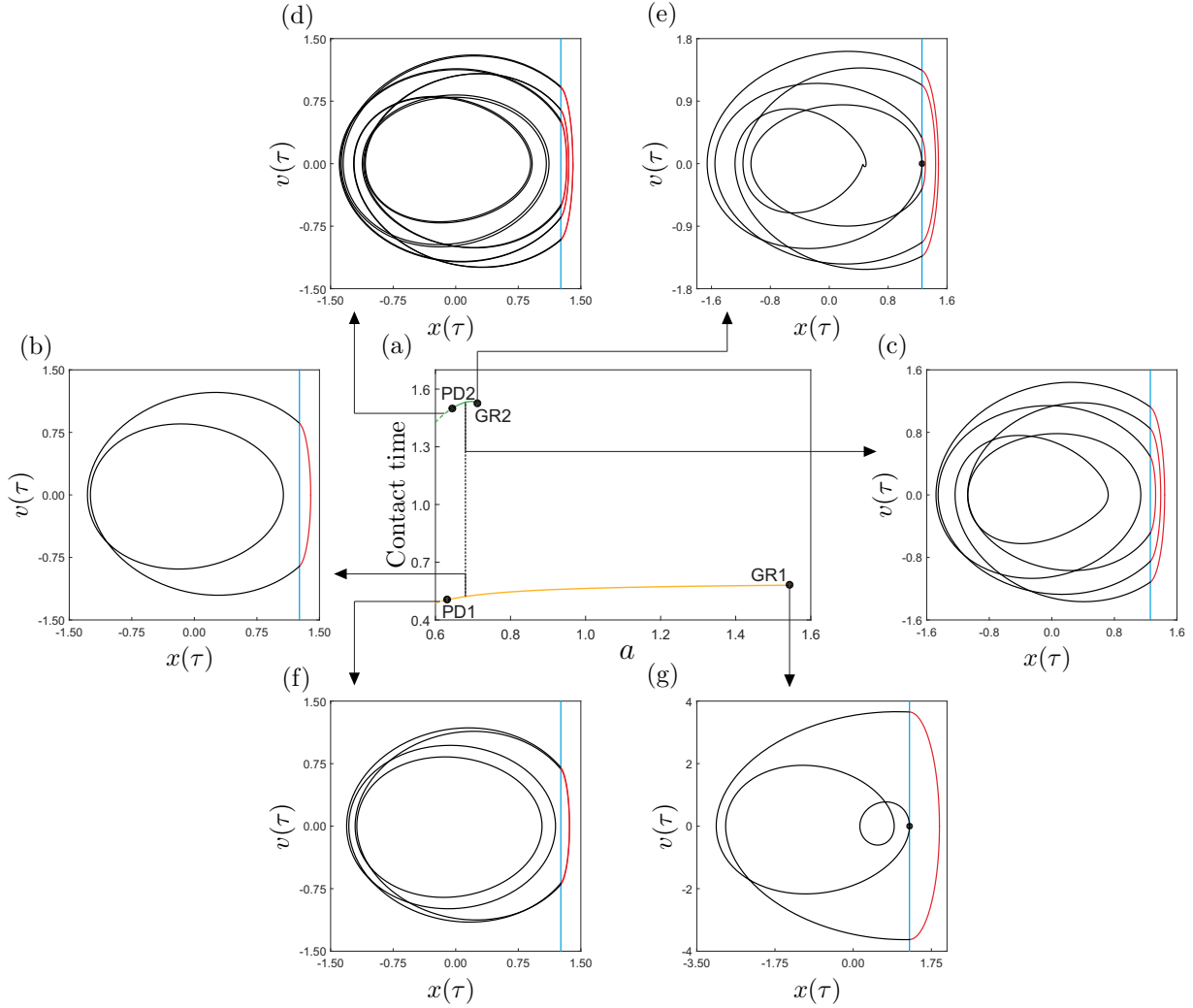}
\caption{(Colour online) (a) One-parameter continuation of the coexisting attractors
shown in panels (b) (yellow branch) and (c) (green branch) with
respect to the excitation amplitude $a$, computed for the parameter
values $\zeta=0.01$, $e=1.26$, $a=0.68$, $\beta=28$ and
$\omega=0.85$ of the impact oscillator \eqref{main-prob}. The
vertical axis shows the \emph{contact time}, the time the impacting
mass stays in contact with the secondary spring per orbital period.
Branches of stable and unstable periodic orbits are depicted with
solid and dashed lines, respectively. The points labeled GR1
($a\approx1.54462$), GR2 ($a\approx0.71258$) and PD1
($a\approx0.63111$), PD2 ($a\approx0.64564$) denote grazing and
period-doubling bifurcations of limit cycles. Panels (g) and (e)
depict periodic solutions corresponding to the grazing bifurcations
points GR1 and GR2, respectively. Here, a dot marks a grazing
contact with the impact boundary $x=e$ (vertical blue line). Panel
(f) presents a period-4 attractor computed for $a=0.63$, while panel
(d) depicts a period-10 solution calculated at $a=0.645$, originated
by the period-doubling bifurcations PD1 and PD2,
respectively.}\label{fig-bif-diag-a}
\end{figure}

As can be seen from the previous section, in order to apply the
proposed control mechanism it is essential to identify parameter
regimes where the considered period-2 and period-5 attractors
maintain both their stability properties and orbit structure. For
this purpose we will carry out a one-parameter continuation of the
underlying attractors with respect to the excitation amplitude $a$,
using the \emph{contact time} as solution measure, i.e.\ the time
the impacting mass stays in contact with the secondary spring per
orbital period.

\begin{figure}[b!]
\centering
\psfrag{X}{\normalsize$x(\tau)$}\psfrag{V}{\normalsize$v(\tau)$}\psfrag{a}{\large$a$}\psfrag{e}{\large$e$}
\psfrag{as}{(a)}\psfrag{b}{(b)}\psfrag{c}{(c)}
\psfrag{d}{(d)}\psfrag{es}{(e)}\psfrag{f}{(f)} \psfrag{g}{(g)}
\includegraphics[width=\textwidth]{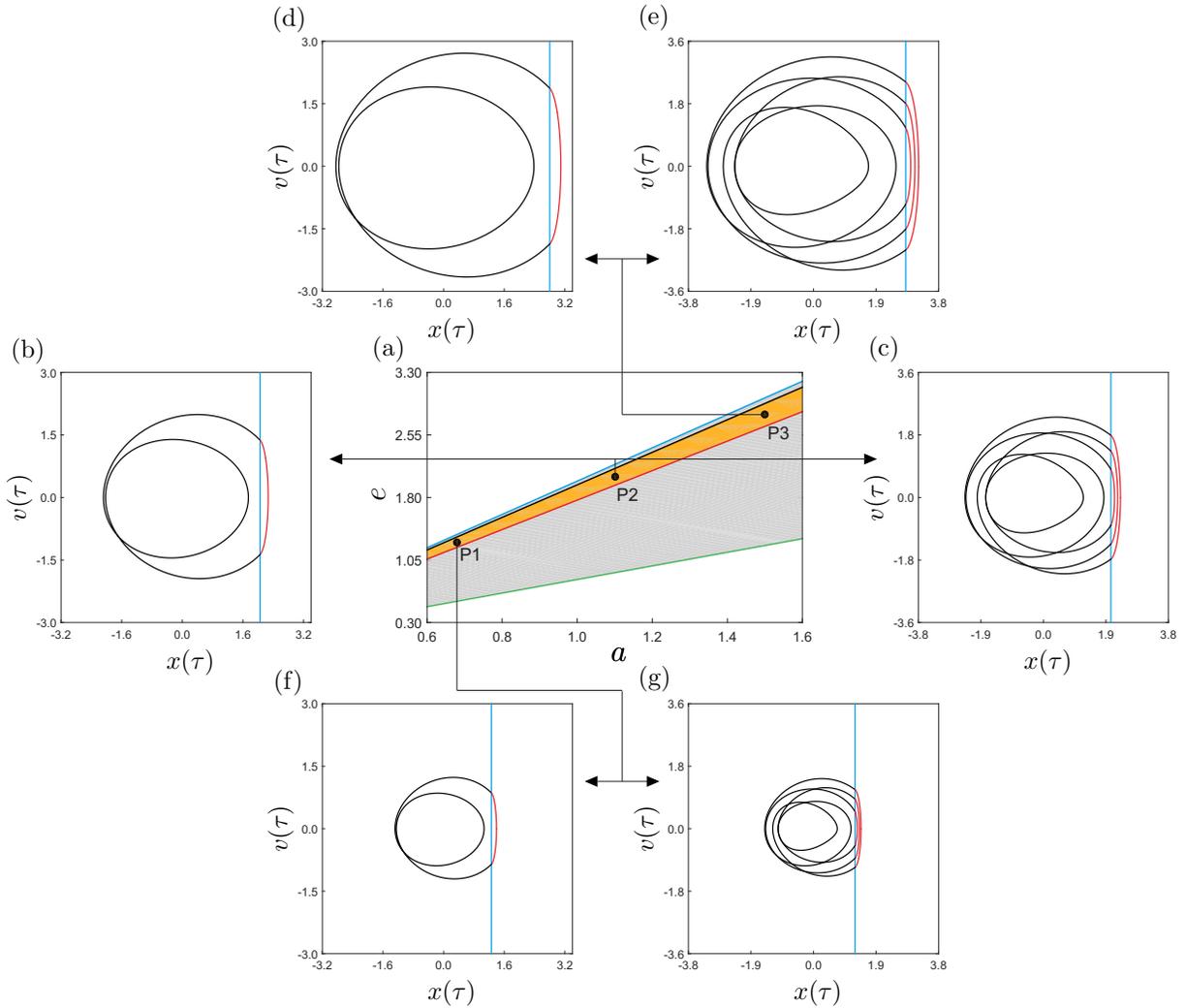}
\caption{(Colour online) (a) Two-parameter continuation of the bifurcation points
PD1 (blue curve), GR1 (green curve), PD2 (black curve) and GR2 (red
curve) found in Fig.\ \ref{fig-bif-diag-a}(a), with respect to the
excitation amplitude $a$ and mass-spring gap $e$. The grey area
represents the parameter region in which the stable period-2
solution of the type shown in panel (b) exists. The yellow region
corresponds to the coexistence of the latter solution type with the
stable period-5 orbit type presented in panel (c). Panels (f)-(g),
(b)-(c) and (d)-(e) represent pairs of coexisting attractors
computed at the test points P1 ($a=0.68$, $e=1.26$), P2 ($a=1.1$,
$e=2.05$) and P3 ($a=1.5$, $e=2.8$).}\label{fig-2par-a-e}
\end{figure}

The result of the process described above is depicted in Fig.\
\ref{fig-bif-diag-a}, which shows yellow and green curves
corresponding to the numerical continuation of the period-2 and
period-5 attractors detected in Fig.\ \ref{coexist}, respectively.
In both cases, the parameter window where the corresponding periodic
solutions remain stable and maintain their orbit structure is
determined by period-doubling and grazing bifurcations of limit
cycles. Specifically, the period-2 solution traced along the yellow
branch loses stability when the excitation amplitude decreases below
$a\approx0.63111$, where the solution undergoes a period-doubling
bifurcation (PD1). Here, the original period-2 orbit (with one
impact per orbital period) becomes unstable and a family of period-4
orbits is born, see for instance the test solution plotted in Fig.\
\ref{fig-bif-diag-a}(f), right after the bifurcation occurs. On the
other hand, when the parameter $a$ increases, a grazing bifurcation
is found at $a\approx1.54462$, where the solution makes tangential
contact with the impact boundary $x=e$, see Fig.\
\ref{fig-bif-diag-a}(g). After this point, a small window of
period-2 solutions with two impacts per orbital period exists, and
they lose stability via a fold bifurcation at $a\approx1.54486$ (not
shown in the diagram). An analogous scenario is found for the
period-5 attractor (with three impacts per orbital period) depicted
in Fig.\ \ref{fig-bif-diag-a}(c). As before, the window of stability
(and orbit structure preservation) for this solution is determined
by the period-doubling bifurcation PD2 ($a\approx0.64564$) and the
grazing point GR2 ($a\approx0.71258$), which results in a
significantly smaller window than the one obtained for the period-2
attractor.

With the results of the one-parameter continuation we are now in
position to determine a parameter region in the $a$-$e$ plane where
the considered period-2 and period-5 attractors maintain both their
stability properties and orbit structure. To this end, we will
perform a two-parameter continuation of the codimension-one
bifurcations detected above. Fig.\ \ref{fig-2par-a-e}(a) shows the
locus of the period-doubling points PD1 (blue curve), PD2 (black
curve) and grazing bifurcations GR1 (green curve), GR2 (red curve)
encountered in Fig.\ \ref{fig-bif-diag-a}(a). In this figure, two
regions are highlighted, in grey and yellow colors. The grey area
represents the parameter region in which the stable period-2
solution (with one impact per orbital period, see panel (b)) exists.
The yellow region corresponds to the coexistence of the latter
solution with the stable period-5 orbit (with three impacts per
orbital period, see panel (c)). Furthermore, several test points
have been selected in order to illustrate the validity of the
highlighted yellow area in the $a$-$e$ plane. Specifically, pairs of
coexisting attractors have been computed at the test points P1
($a=0.68$, $e=1.26$), P2 ($a=1.1$, $e=2.05$) and P3 ($a=1.5$,
$e=2.8$), see panels (b)--(g) in Fig.\ \ref{fig-2par-a-e}. In this
way, the yellow area can be used as a reference for the
applicability of the proposed control scheme, so as to guarantee
that the parametric perturbations do not bring the system to a
regime where either of the considered attractors lose stability or
the intended orbit structure.

\subsection{Controlling multiple coexisting attractors}\label{sec-three-coex}

\begin{figure}[h!]
\centering
\includegraphics[width=0.85\textwidth]{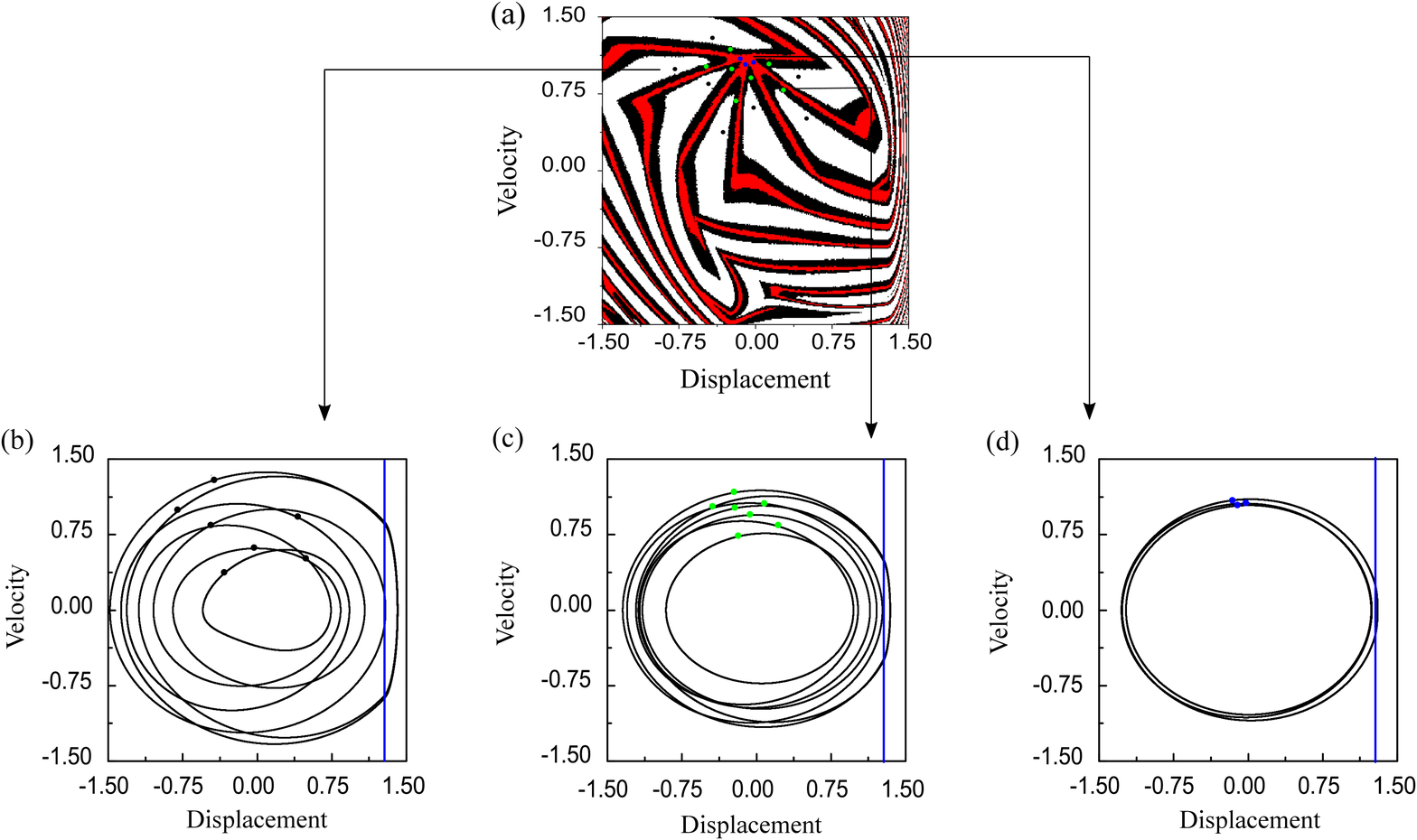}
\caption{(Colour online) (a) Basins of attraction of the impacting
system computed for $\zeta= 0.01, ~e = 1.28, ~a = 0.49, ~\beta = 28 ~
\textrm{and}~ \omega= 0.8528$. Black dots, green dots and blue dots denote the period-$7$
attractor with large amplitude and white basin, the period-$7$
attractor with small amplitude and black basin, and the period-$3$ attractor with red basin. (b), (c) and (d) present the trajectories
of the period-$7$ with large amplitude, the period-$7$ with small amplitude, and the period-$3$ attractors on the phase plane,
where blue lines indicate the impact boundary.}\label{three_coexist}
\end{figure}

\begin{figure}[h!]
\centering
\includegraphics[width=0.9\textwidth]{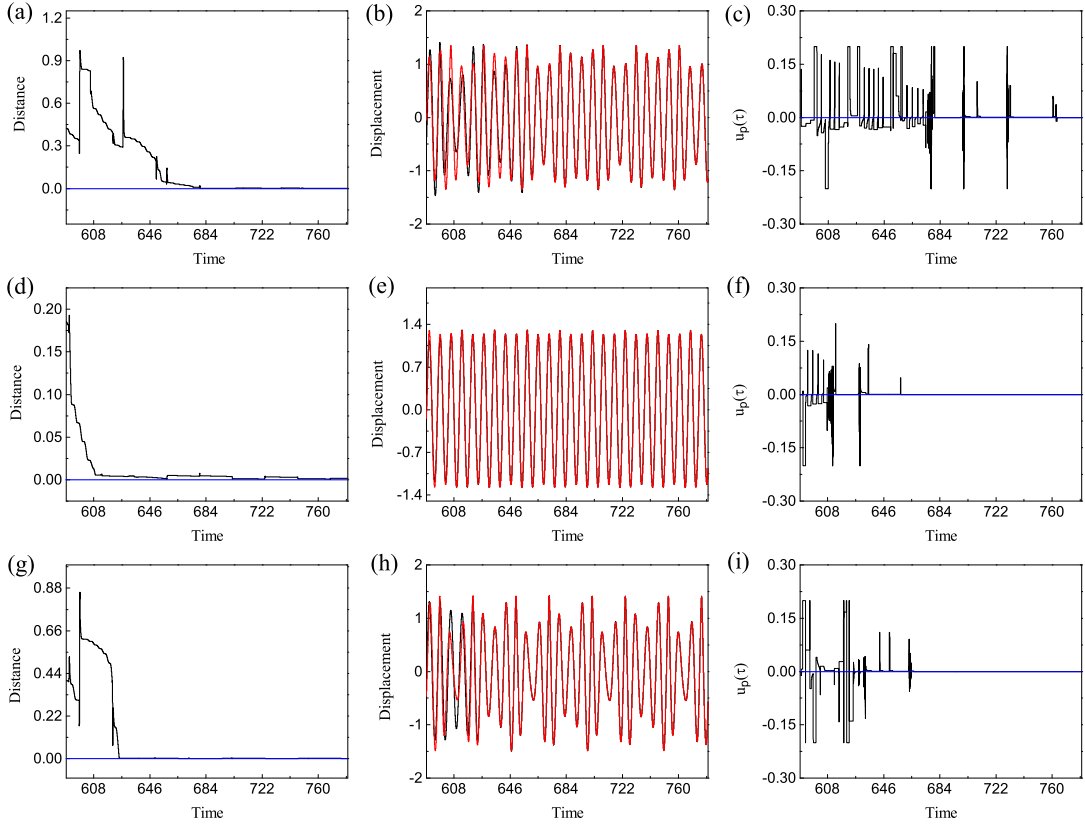}
\caption{(Colour online) (a)-(c) present the controlled details from the initial attractor (period-$7$ with large amplitude) to the desired attractor (period-$7$ attractor with small amplitude) of the system under the nonlinear control strategy (Algorithm 2) by varying the amplitude of excitation with $M_{a,1}=0.2$ and $M_{a,2}=10$.
(d)-(f) present the controlled details from the initial attractor (period-$7$ with small amplitude) to the desired attractor (period-$3$ attractor) of the system under the same nonlinear control strategy and conditions.
(g)-(i) present the controlled details from the initial attractor (period-$3$ attractor) to the desired attractor (period-$7$ attractor with large amplitude) of the system under the same nonlinear control strategy and conditions.
(a), (d) and (g) show the time histories of the distance between the desired and the controlled trajectories in 2-norm.
(b), (e) and (h) show the time histories of the desired (red line) and the current (black line) displacements.
(c), (f) and (i) show the time histories of the control sequence generated by the nonlinear control strategy.
Blue lines mark the zero reference.
} \label{three_attractor}
\end{figure}

Here, we consider achieving the switch among three coexisting attractors, as shown in Fig.\ \ref{three_coexist}, by the nonlinear control strategy (Algorithm 2) through varying the amplitude. At the beginning, the control of the impacting system from the initial attractor (period-$7$ attractor with large amplitude) to desired attractor (period-$7$ attractor with small amplitude) is considered, as shown in Fig.\ \ref{three_attractor}(a)-(c).
In details, in Fig.\ \ref{three_attractor}(a), the distance between the desired and the controlled trajectories in $2$-norm was decreased and finally was reduced to nil.
Fig.\ \ref{three_attractor}(c) presents the time history of the control signal $u_{p}(\tau)$ that initiated from $\tau=589.417$ and terminated at $\tau=763.2234$.
The transition on the displacement can be observed from Fig.\ \ref{three_attractor}(b).
Secondly, the control from the period-$7$ attractor with small amplitude to the period-$3$ attractor by varying the amplitude of excitation $a$ is presented in Fig.\ \ref{three_attractor}(d)-(f).
Specifically, the control strategy was applied at $\tau=589.417$ and was switched off at $\tau=657.5665$.
During this period,  the trajectory experienced a transition, as shown in Fig.\ \ref{three_attractor}(e), and the distance between the controlled and the desired trajectories was decreased from $0.182$ to $0.001$, as shown in Fig.\ \ref{three_attractor}(d).
In the meanwhile, the control signal did not change significantly and reduced to nil gradually.
After that, the controlled trajectory approached the desired trajectory spontaneously.
Finally, Fig.\ \ref{three_attractor}(g)-(i) present the control result from the period-$3$ attractor to period-$7$ attractor with large amplitude by varying the amplitude.
The control strategy was also applied at $\tau=589.417$, and switched off at $\tau=666.3881$ as shown in Fig.\ \ref{three_attractor}(i).
Within the same time duration, the trajectory of the controlled system approached to the desired trajectory in Fig.\ \ref{three_attractor}(h), and the distance between the controlled and the desired trajectories was decreased from $0.403$ to $0.001$ as shown in Fig.\ \ref{three_attractor}(g).
The above processes show that the proposed method can present good performances on controlling more coexisting attractors.

\section{Control of smooth dynamical systems}\label{sec-smooth}
\subsection{Nonlinear control}
\begin{figure}[h!]
\centering
\includegraphics[width=0.55\textwidth]{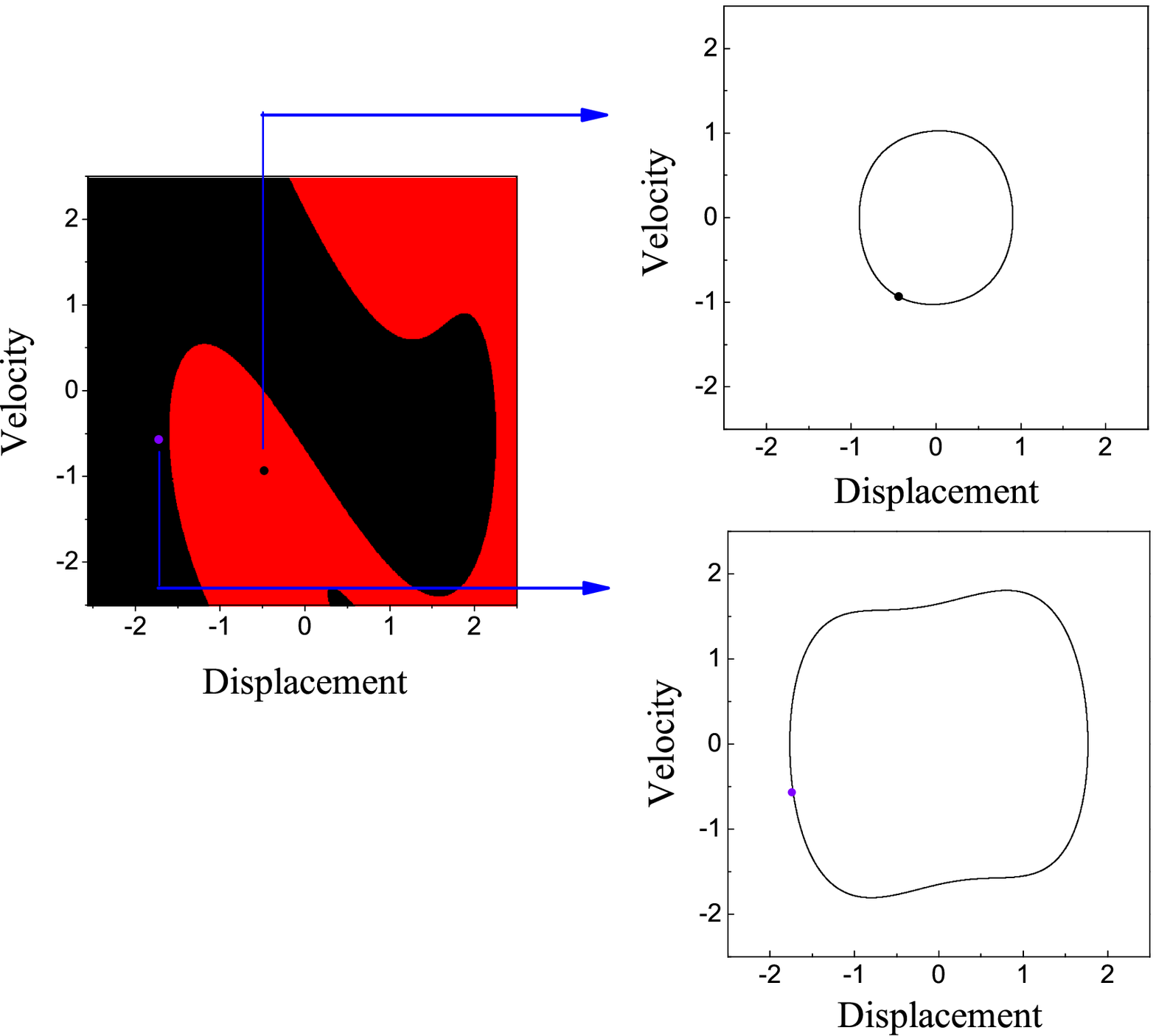}
\caption{(Colour online) Basins of attraction of the Duffing system
computed for $\Gamma = 1.9$, $\omega = 1.2$, $p_{1} = 0.9$ and
$p_{2}=1$. Black dot denotes the period-$1$ small amplitude
attractor with red basin, and violate dot represents the period-$1$
large amplitude attractor with black basin. The right panels present
the trajectories of the two period-$1$ attractors on the phase
plane.} \label{coexistduffing}
\end{figure}

In this section, the Duffing oscillator representing smooth dynamical
systems is employed to test the versatility of the proposed control
method. The Duffing system, which is known to have many coexisting
attractors without control, can be described by
\begin{equation}\label{duffing-prob}
\begin{cases}
\dot{Y}(\tau)=F_\mathrm{du}(\tau,Y(\tau),u_\mathrm{p}(\tau)), \\
Y_{0}=Y(\tau_{0}),
\end{cases}
\end{equation}
where $Y(\tau):=(x(\tau),v(\tau))^{T}$, and
\begin{align*}
F_\mathrm{du}(\tau,Y,u_\mathrm{p}):=\left[\begin{array}{cc}
  0 \\
\Gamma \sin(\omega\tau) \\
\end{array}\right]+
\left[\begin{array}{cc}
0 & 1 \\
1 & -p_{1} \\
\end{array}\right]Y
+Y^{T}(\tau)\left[\begin{array}{cc}
  1\\
0 \\
\end{array}\right]Y^{T}\left[\begin{array}{cc}
  1\\
0 \\
\end{array}\right]
\left[\begin{array}{cc}
  0 & 0\\
-(p_{2}+u_\mathrm{p}) & 0  \\
\end{array}\right]Y.
\end{align*}
The following parameters: $\Gamma = 1.9$, $\omega = 1.2$, $p_{1} =
0.9$ and $p_{2}=1$ were considered in this study. At these parameter
values, the system without control ($u_\mathrm{p}=0$) has two coexisting attractors, depicted in
Fig.~\ref{coexistduffing}, which are a period-1 small and a large
amplitude attractors with their Poincar\'{e} sections denoted by
black and violet dots, respectively.

\begin{figure}[h!]
\centering
\includegraphics[width=0.82\textwidth]{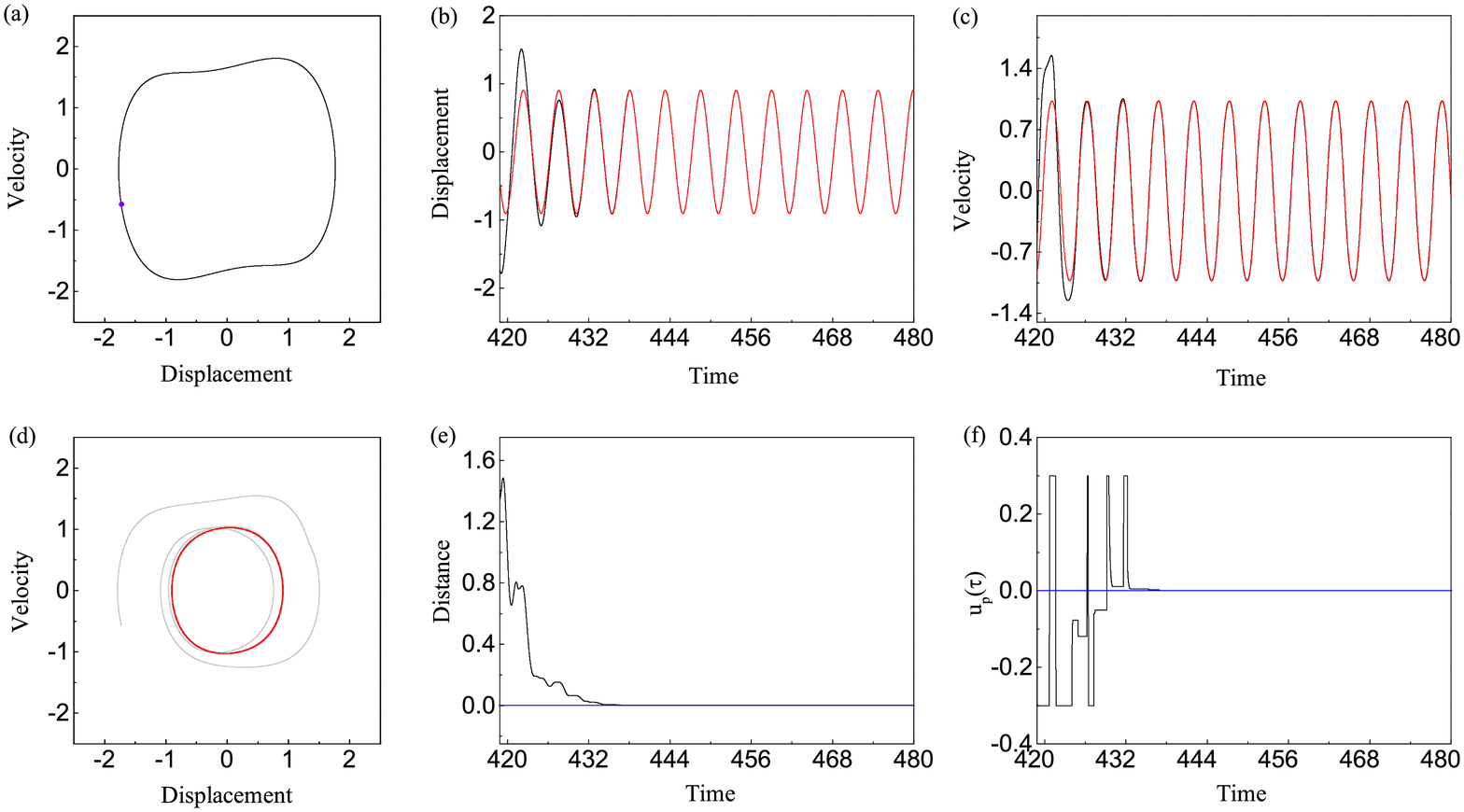}
\caption{(Colour online) (a) The large amplitude period-1 response
on the phase plane with the Poincar\'{e} section denoted by violate
dot. (b) Time histories of the desired (red line) and the current
(black line) displacements of the system under the nonlinear
control strategy (Algorithm 2) by varying the stiffness of the
nonlinear spring with $M_\mathrm{p_{2},1}=0.3$ and
$M_\mathrm{p_{2},2}=10$. (c) Time histories of the desired (red
line) and the current (black line) velocities of the system. (d)
Trajectory of the system on the phase plane under the
nonlinear control strategy, where grey and red lines
represent the transient and the steady-state responses,
respectively. (e) Time history of the distance between the desired
and the controlled trajectories in $2$-norm. (f) Time history of the
control sequence generated by the nonlinear control
strategy. Blue lines in (a) and (d) indicate the impact boundary,
while the blue lines in (e) and (f) mark the zero reference. The
result was computed for $\Gamma = 1.9$, $\omega = 1.2$, $p_{1} =
0.9$ and $p_{2}=1$.} \label{four1}
\end{figure}

The control aims for system (\ref{duffing-prob}) to switch the
two stable attractors shown in Fig.~\ref{coexistduffing} by varying
the stiffness of the nonlinear spring $p_{2}$. The control result for the switching from the large to the small amplitude attractor is shown in Fig.~\ref{four1}, where the control
strategy was applied to the original attractor at $\tau=418.879$,
and the controlled trajectory experienced a transition until
$\tau=439$. During this time the 2-norm distance between the
control and the desired trajectories was reduced from $1.338$ to
$0.001$. The control signal reached the maximum
$u_\mathrm{p}(\tau)=0.3$ and the minimum $u_\mathrm{p}(\tau)=-0.3$
for several times before it was switched off.

\begin{figure}[h!]
\centering
\includegraphics[width=0.82\textwidth]{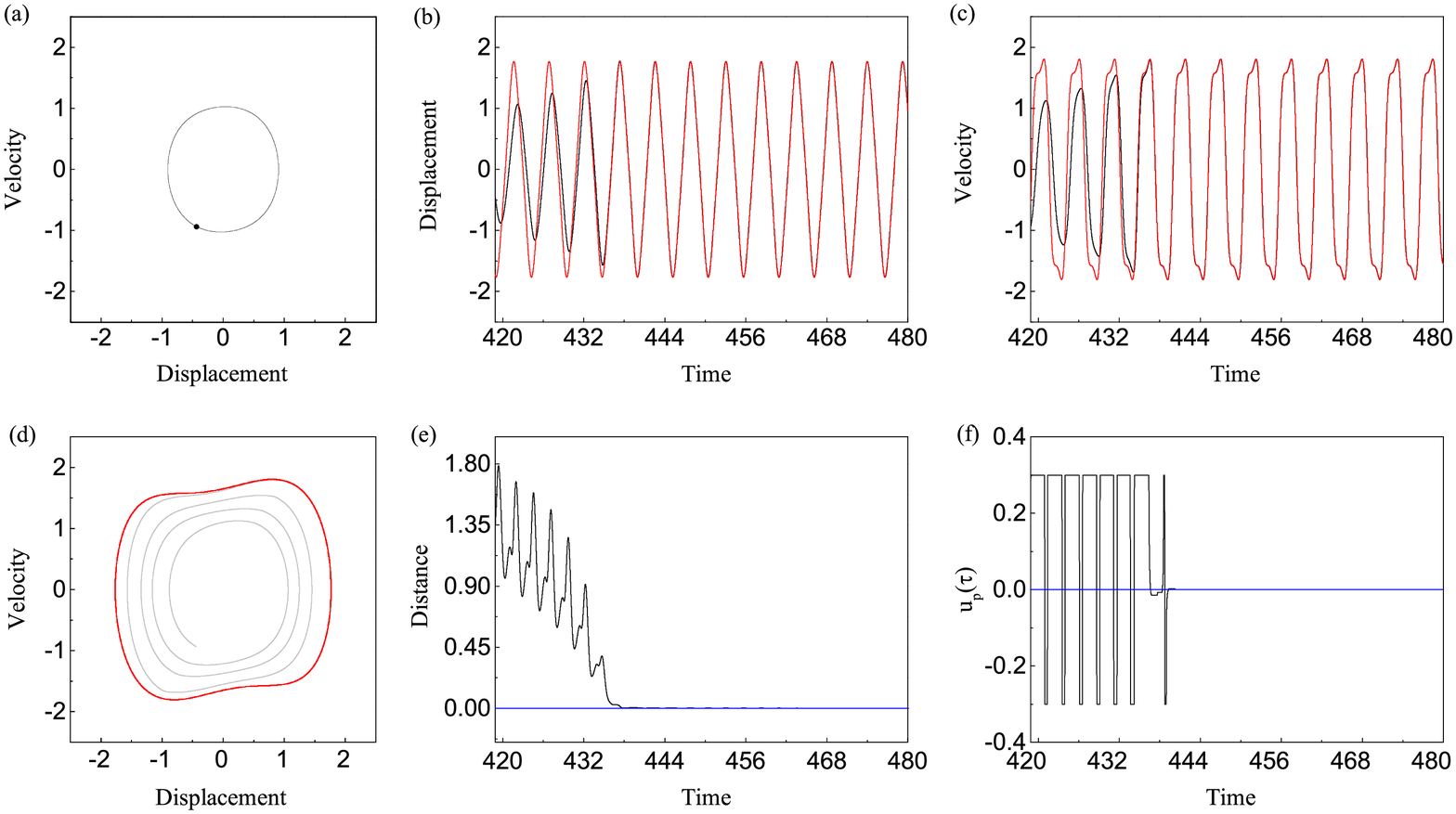}
\caption{(Colour online) (a) The small amplitude period-1 response
on the phase plane with the Poincar\'{e} section denoted by black
dot. (b) Time histories of the desired (red line) and the current
(black line) displacements of the system under the nonlinear
control strategy (Algorithm 2) by varying the stiffness of the
nonlinear spring with $M_\mathrm{p_{2},1}=0.3$ and
$M_\mathrm{p_{2},2}=10$. (c) Time histories of the desired (red
line) and the current (black line) velocities of the system. (d)
Trajectory of the system on the phase plane under the
nonlinear control strategy, where grey and red lines
represent the transient and the steady-state responses,
respectively. (e) Time history of the distance between the desired
and the controlled trajectories in $2$-norm. (f) Time history of the
control sequence generated by the nonlinear control
strategy. Blue lines in (a) and (d) indicate the impact boundary,
while the blue lines in (e) and (f) mark the zero reference. The
result was computed for $\Gamma = 1.9$, $\omega = 1.2$, $p_{1} =
0.9$ and $p_{2}=1$.} \label{fourdifferent1}
\end{figure}

The control from the small to the large amplitude attractor is
presented in Fig.~\ref{fourdifferent1}, where the control strategy
was switched on at $\tau=418.879$ and was switched off at
$\tau=440.246$ when the distance between the two trajectories was
decreased to $0.001$. Compared to the switching in Fig.~\ref{four1},
the transition from the small to the large amplitude attractor took
a longer time.

\subsection{Bifurcation analysis of the coexisting attractors}

Analogous to Section \ref{sec-coex-IO}, in this section our main
concern will be to study in detail the effect of the control
parameters $p_{1}$, $p_{2}$ on the small- and high-amplitude
oscillations of the Duffing system \eqref{duffing-prob}, see Fig.\
\ref{coexistduffing}. To this end, we will employ path-following
methods for limit cycles, implemented via the continuation platform
COCO \cite{dankowicz2013recipes}, along with its routines for
bifurcation detection and two-parameter continuation of
codimension-1 bifurcations.

\begin{figure}[h!]
\centering
\psfrag{X}{\normalsize$x(\tau)$}\psfrag{V}{\normalsize$v(\tau)$}\psfrag{p1}{\normalsize$p_{1}$}\psfrag{p2}{\normalsize$p_{2}$}
\psfrag{P2P}{\normalsize$A_{\mbox{\tiny P2P}}$}
\psfrag{as}{(a)}\psfrag{b}{(b)}\psfrag{c}{(c)} \psfrag{d}{(d)}
\includegraphics[width=\textwidth]{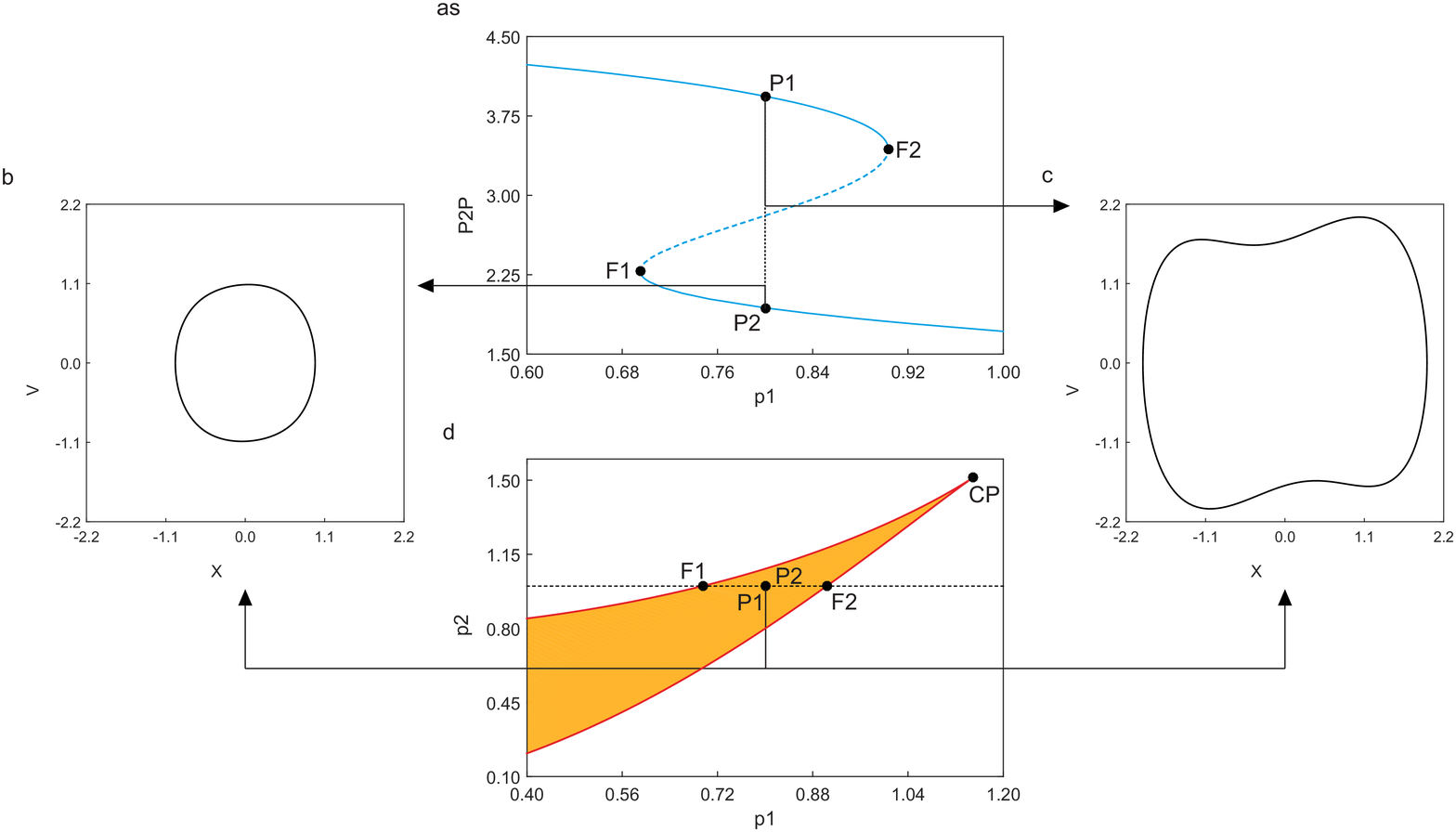}
\caption{(Colour online) (a) One-parameter continuation of the periodic response of
the Duffing oscillator \eqref{duffing-prob} with respect to $p_{1}$,
computed for the parameter values $\Gamma=1.9$, $\omega=1.2$ and
$p_{2}=1$. The vertical axis shows the peak-to-peak amplitude
$A_{\mbox{\tiny P2P}}$ of the $x$-component. Branches of stable and
unstable periodic orbits are depicted with solid and dashed lines,
respectively. The points labeled F1 ($p_{1}\approx0.69494$) and F2
($p_{1}\approx0.90352$) stand for fold bifurcations of limit cycles.
Panels (b) and (c) depict stable coexisting solutions computed at
the test points P1 and P2 ($p_{1}=0.8$), respectively, shown in
panel (a). (d) Two-parameter continuation of the fold points found
in panel (a), with respect to $p_{1}$ and $p_{2}$. Here, the label
CP represents a cusp bifurcation ($p_{1}\approx1.14902$,
$p_{2}\approx1.51194$). The intersections of the horizontal dashed
line ($p_{2}=1$) with the bifurcation diagram correspond to the fold
bifurcations F1 and F2 shown in panel (a). The yellow area
represents the parameter region in which the stable periodic
solutions of the type shown in panels (b) and (c) coexist.
}\label{fig-bif-duff}
\end{figure}

The starting point for our study is the high-amplitude periodic
solution shown in Fig.\ \ref{fig-bif-duff}(c), computed for
$p_{1}=0.8$. Panel (a) presents the result of the numerical
continuation of this orbit with respect to the control parameter
$p_{1}$. In this diagram, changes of stability are detected, which
are marked with solid (for stable solutions) and dashed (unstable
solutions) lines. The window of stability of the the high-amplitude
orbit is bounded from above by the fold bifurcation F2
($p_{1}\approx0.90352$). At this point, a branch of unstable
periodic solutions is born, which finishes at the fold point F1
($p_{1}\approx0.69494$). Here, a family of stable oscillations
emerges, corresponding to small-amplitude periodic orbits as can be
seen at the test point P2, see Fig.\ \ref{fig-bif-duff}(b).
Consequently, the bifurcation points F1 and F2 defines a parameter
window where both attractors coexist.

Next, we will carry out a two-parameter continuation of the fold
points detected above in order to determine a region in the
$p_{1}$-$p_{2}$ plane where the small- and high-amplitude attractors
of the Duffing system coexist. The result of this numerical process
is presented in Fig.\ \ref{fig-bif-duff}(d), where the red curve
stands for a locus of fold bifurcations of limit cycles. In this
picture, the yellow area enclosed by the fold curve represents the
parameter regime where small- and high-amplitude oscillations
coexist. The intersections of the horizontal dashed line ($p_{2}=1$)
with the bifurcation diagram correspond to the fold bifurcations F1
and F2 found in Fig.\ \ref{fig-bif-duff}(a). Furthermore, the
numerical computations reveal the presence of a codimension-2 point
$(p_{1},p_{2})\approx(1.14902,1.51194)$ (CP), where two branches of
fold points (those corresponding to F1 and F2) join together via a
cusp singularity. In this way, it is possible to determine
boundaries in the considered parameter region for the application of
the control mechanism proposed in this work.

\section{Concluding remarks}\label{sec-conclusion}

This paper studied a new control method for switching stable coexisting attractors of
non-autonomous smooth and non-smooth dynamical systems. Our control aim was to control an undesired coexisting attractor to a desired one by modulating a system parameter without affecting the original property of the system. To examine the proposed control concept, we implemented two control strategies with finite sampling step, namely the linear and nonlinear control strategies, where one was implemented through the external control input and the other one was applied via a system parameter. In the first part of our simulation work, two multistable scenarios (one is coexisting two attractors and another is coexisting three attractors) of the impact oscillator were studied. Our simulations show that both control strategies are effective for switching the stable coexisting attractors in the impact oscillator. The effective control region for the control parameters (excitation amplitude $a$ and mass-spring gap $e$) was also found by employing the path-following methods for piecewise-smooth dynamical systems. In the second part of the simulation work, we implemented the nonlinear control strategy to the Duffing oscillator for switching a period-1 small and a period-1 large amplitude attractors. Path-following methods for limit cycles were used to identify the effective control region of the proposed control method.

Compared with the classical delay feedback control proposed by Pyragas \cite{pyragas1992continuous}, our proposed methods have many advantages on the control of coexisting attractors.
First of all, the main advantage of the proposed method for control of coexisting attractors compared to delayed feedback control is its transient behaviour.
Delayed feedback control \emph{locally} stabilizes a target periodic orbit that is unstable without control. Its behaviour for initial conditions far away from the target is not considered in its design and, thus, can lead to undesirable arbitrarily long transients. In contrast, the proposed method contains terms explicitly driving the trajectory toward its target.
Secondly, especially for periodically forced systems the delay feedback control is infeasible for switching from an initial attractor to a target with equal (or multiple) period. This is apparent in the examples in sections 3.4 (period 7 to period 7), 4.1 and 4.2 (both period 1 to period 1), where delayed feedback would be non-invasive (and, hence, ineffective) on the initial condition.
Thirdly, the proposed methods do not result in a significant change on the system's states due to restricting the varying rate of the control signal.
On the contrary, when the delay feedback control is just introduced into the system, the system's states have to witness a significant change due to the value of the control is very large at the beginning.
Finally, the nonlinear control strategy depends only on the original properties of system parameter, and does not need to introduce any external inputs.
Besides the above differences, the common point between the proposed methods and delay feedback control is needing the information of derivatives of the right-hand sides of the ODE.
In details, although the proposed methods require estimates of these derivatives to generate the control signal in its implementation and the delay feedback control does not, the feasible range of control gains for which delayed feedback control is stabilizing is still limited by the same derivatives, due to the stability of controlled trajectory depends on the Jacobian \cite{zhang2021controlling}.


Future works include optimisation, experimental validation and implementation of the proposed control method in more complex multistable scenarios in both smooth and non-smooth dynamical systems.

\section*{Acknowledgements}
This work has been supported by EPSRC under Grant No. EP/P023983/1.
Mr Zhi Zhang would like to acknowledge the financial support from
the University of Exeter for his Exeter International Excellence
Scholarship. Prof. Jan Sieber's research is supported by EPSRC
Fellowship EP/N023544/1 and EPSRC grant EP/V04687X/1.

\section*{Compliance with ethical standards}

\section*{Conflict of interest}
The authors declare that they have no conflict of interest
concerning the publication of this manuscript.

\section*{Data accessibility}
The datasets generated and analysed during the current study are
available from the corresponding author on reasonable request.

\section*{Appendix}
\textbf{Proof of Theorem \ref{Main_1}}:
Assume that there exists a sufficiently small $h>0$, the interval
$[\tau_{0},\tau^{*}]=\bigcup_{i=0}^{n^{*}}[\tau_{i},\tau_{i+1}]$,
where $\tau_{i+1}=\tau_{i}+h$ and $i=1,\cdots,n^{*}$.
For $\tau\in[\tau_{i},\tau_{i}+h]$, the distance vector
can be defined as
\begin{displaymath}
  d(\tau):=Y_{\mathrm{d}}(\tau_{0})+\int_{\tau_{0}}^{\tau}F(\tau_\mathrm{s},
  Y_{\mathrm{d}}(\tau_\mathrm{s})) \mathrm{d}\tau_\mathrm{s}-Y_{\mathrm{c}}(\tau_{0})-
  \int_{\tau_{0}}^{\tau}(F(\tau_\mathrm{s},Y_{\mathrm{u}}(\tau_\mathrm{s})) +
  U(\tau_\mathrm{s}))\mathrm{d}\tau_\mathrm{s}.
\end{displaymath}
Thus,
\begin{align}\label{distan-1}
\langle d(\tau), d(\tau) \rangle=&\langle Y_{\mathrm{d}}(\tau_\mathrm{i})+
\int_{\tau_\mathrm{i}}^{\tau}F(\tau_\mathrm{s},Y_{\mathrm{d}}(\tau_\mathrm{s}))
 \mathrm{d} \tau_\mathrm{s}-Y_{\mathrm{u}}(\tau_\mathrm{i})-\int_{\tau_\mathrm{i}}^{\tau}
 (F(\tau_\mathrm{s},Y_{\mathrm{u}}(\tau_\mathrm{s})) +U(\tau_\mathrm{s}))\mathrm{d} \tau_\mathrm{s},\nonumber\\
&
Y_{\mathrm{d}}(\tau_\mathrm{i})+\int_{\tau_\mathrm{i}}^{\tau}F(\tau_\mathrm{s},
Y_{\mathrm{d}}(\tau_\mathrm{s}))
\mathrm{d}
\tau_\mathrm{s}-Y_{\mathrm{u}}(\tau_\mathrm{i})-\int_{\tau_\mathrm{i}}^{\tau}(F(\tau_\mathrm{s},Y_{\mathrm{u}}(\tau_\mathrm{s}))
+U(\tau_\mathrm{s}))\mathrm{d} \tau_\mathrm{s}    \rangle,
\end{align}
where
$Y_{\mathrm{u}}(\tau_\mathrm{i}):=Y_{\mathrm{c}}(\tau_{0})+\int_{\tau_{0}}^{\tau_\mathrm{i}}(F(\tau_\mathrm{s},Y_{\mathrm{u}}(\tau_\mathrm{s}))
+U(\tau_\mathrm{s}))\mathrm{d}\tau_\mathrm{s}$. Since
$U(\tau)=U(\tau_\mathrm{i})+\dot{U}(\tau_\mathrm{i})(\tau-\tau_\mathrm{i})$,
where $\tau\in [\tau_\mathrm{i},\tau_\mathrm{i}+h]$, and
Eq.~\eqref{dist:form}, it gives
\begin{align*}
\langle d(\tau), d(\tau) \rangle=&\, \langle Y_{\mathrm{d}}(\tau_\mathrm{i})-Y_{\mathrm{u}}(\tau_\mathrm{i})+
\int_{\tau_\mathrm{i}}^{\tau}(F(\tau_\mathrm{s},Y_{\mathrm{d}}(\tau_\mathrm{s})) -
F(\tau_\mathrm{s},Y_{\mathrm{u}}(\tau_\mathrm{s})) -U(\tau_\mathrm{s}))\mathrm{d} \tau_\mathrm{s},\nonumber\\
& \, Y_{\mathrm{d}}(\tau_\mathrm{i})-Y_{\mathrm{u}}(\tau_\mathrm{i})+\int_{\tau_\mathrm{i}}^{\tau}
(F(\tau_\mathrm{s},Y_{\mathrm{d}}(\tau_\mathrm{s})) -F(\tau_\mathrm{s},Y_{\mathrm{u}}(\tau_\mathrm{s})) -
U(\tau_\mathrm{s}))\mathrm{d} \tau_\mathrm{s}    \rangle\\
=&\, \langle d(\tau_\mathrm{i}), d(\tau_\mathrm{i}) \rangle + 2\langle d(\tau_\mathrm{i}),
 \int_{\tau_\mathrm{i}}^{\tau}(F(\tau_\mathrm{s},Y_{\mathrm{d}}(\tau_\mathrm{s}))-
 F(\tau_\mathrm{s},Y_{\mathrm{u}}(\tau_\mathrm{s}))-U(\tau_\mathrm{i})) \mathrm{d} \tau_\mathrm{s} \rangle\\
&-2\langle d(\tau_\mathrm{i}), \int_{\tau_\mathrm{i}}^{\tau} \dot{U}(\tau_\mathrm{i})
(\tau_\mathrm{s}-\tau_\mathrm{i})\mathrm{d} \tau_\mathrm{s} \rangle\\
&+\langle \int_{\tau_\mathrm{i}}^{\tau}(F(\tau_\mathrm{s},Y_{\mathrm{d}}(\tau_\mathrm{s})) -
F(\tau_\mathrm{s},Y_{\mathrm{u}}(\tau_\mathrm{s}))-U(\tau_\mathrm{i})) \mathrm{d} \tau_\mathrm{s},
\int_{\tau_\mathrm{i}}^{\tau}(F(\tau_\mathrm{s},Y_{\mathrm{d}}(\tau_\mathrm{s})) -
F(\tau_\mathrm{s},Y_{\mathrm{u}}(\tau_\mathrm{s}))-U(\tau_\mathrm{i})) \mathrm{d} \tau_\mathrm{s} \rangle\\
&-2\langle \int_{\tau_\mathrm{i}}^{\tau}(F(\tau_\mathrm{s},Y_{\mathrm{d}}(\tau_\mathrm{s})) -
F(\tau_\mathrm{s},Y_{\mathrm{u}}(\tau_\mathrm{s}))-U(\tau_\mathrm{i})) \mathrm{d} \tau_\mathrm{s},
 \int_{\tau_\mathrm{i}}^{\tau} \dot{U}(\tau_\mathrm{i})(\tau_\mathrm{s}-\tau_\mathrm{i})\mathrm{d}
  \tau_\mathrm{s} \rangle\\
&+\langle\int_{\tau_\mathrm{i}}^{\tau}
\dot{U}(\tau_\mathrm{i})(\tau_\mathrm{s}-\tau_\mathrm{i})\mathrm{d}
\tau_\mathrm{s}, \int_{\tau_\mathrm{i}}^{\tau}
\dot{U}(\tau_\mathrm{i})(\tau_\mathrm{s}-\tau_\mathrm{i})\mathrm{d}
\tau_\mathrm{s}\rangle.
\end{align*}

When $\tau=\tau_{i+1}$, it can obtain that,
\begin{align}\label{onestep_condition}
&\langle d(\tau_{i+1}), d(\tau_{i+1}) \rangle-\langle d(\tau_\mathrm{i}), d(\tau_\mathrm{i})
 \rangle \nonumber\\
&\qquad\qquad=2\langle d(\tau_\mathrm{i}),
F(\tau_\mathrm{i},Y_{\mathrm{d}}(\tau_\mathrm{i}))-F(\tau_\mathrm{i},Y_{\mathrm{u}}(\tau_\mathrm{i}))-U(\tau_\mathrm{i})
\rangle h
-\langle d(\tau_\mathrm{i}),  \dot{U}(\tau_\mathrm{i}) \rangle h^{2} \nonumber\\
&\qquad\qquad\quad+\langle F(\tau_\mathrm{i},Y_{\mathrm{d}}(\tau_\mathrm{i})) -
F(\tau_\mathrm{i},Y_{\mathrm{u}}(\tau_\mathrm{i}))-U(\tau_\mathrm{i}) ,  F(\tau_\mathrm{i},Y_{\mathrm{d}}
(\tau_\mathrm{i})) -F(\tau_\mathrm{i},Y_{\mathrm{u}}(\tau_\mathrm{i}))-U(\tau_\mathrm{i}) \rangle h^{2} \nonumber\\
&\qquad\qquad\quad-\langle
F(\tau_\mathrm{i},Y_{\mathrm{d}}(\tau_\mathrm{i}))
-F(\tau_\mathrm{i},Y(\tau_\mathrm{i}))-U(\tau_\mathrm{i}) ,
\dot{U}(\tau_\mathrm{i}) \rangle h^{3}
+\tfrac{1}{4}\langle  \dot{U}(\tau_\mathrm{i}), \dot{U}(\tau_\mathrm{i}) \rangle h^{4} \nonumber\\
&\qquad\qquad=2\langle d(\tau_\mathrm{i}),
F(\tau_\mathrm{i},Y_{\mathrm{d}}(\tau_\mathrm{i}))-F(\tau_\mathrm{i},Y_{\mathrm{u}}(\tau_\mathrm{i}))-U(\tau_\mathrm{i})
\rangle h
-\langle d(\tau_\mathrm{i}),  \dot{U}(\tau_\mathrm{i}) \rangle h^{2} \\
&\qquad\qquad\quad+\langle
F(\tau_\mathrm{i},Y_{\mathrm{d}}(\tau_\mathrm{i}))
-F(\tau_\mathrm{i},Y_{\mathrm{u}}(\tau_\mathrm{i}))-U(\tau_\mathrm{i}) ,
F(\tau_\mathrm{i},Y_{\mathrm{d}}(\tau_\mathrm{i}))
-F(\tau_\mathrm{i},Y_{\mathrm{u}}(\tau_\mathrm{i}))-U(\tau_\mathrm{i}) \rangle
h^{2}+O(h^{3}). \nonumber
\end{align}
If $\langle d(\tau_{i+1}), d(\tau_{i+1}) \rangle-\langle
d(\tau_\mathrm{i}), d(\tau_\mathrm{i}) \rangle$ is not positive, the
inequality (\ref{condition}) in \textbf {Step 2} is obtained. By
repeating $n^{*}$ times, we can obtain
\begin{align*}
&\langle d(\tau_{n^{*}}), d(\tau_{n^{*}}) \rangle-\langle d(\tau_{0}),d(\tau_{0}) \rangle\\
&\qquad\qquad=\sum_{i=0}^{n^{*}-1}\big[\langle d(\tau_{i+1}), d(\tau_{i+1}) \rangle-\langle
d(\tau_\mathrm{i}), d(\tau_\mathrm{i}) \rangle\big] \\
&\qquad\qquad=\sum_{i=0}^{n^{*}}\big[2\langle d(\tau_\mathrm{i}),
F(\tau_\mathrm{i},
Y_{\mathrm{d}}(\tau_\mathrm{i}))-F(\tau_\mathrm{i},
Y_{\mathrm{u}}(\tau_\mathrm{i})) -U(\tau_\mathrm{i})\rangle h-
\langle d(\tau_\mathrm{i}), \dot{U}(\tau_\mathrm{i}) \rangle h^2\\
&\qquad\qquad\quad+\langle F(\tau_\mathrm{i},
Y_{\mathrm{d}}(\tau_\mathrm{i}))-F(\tau_\mathrm{i},
Y_{\mathrm{u}}(\tau_\mathrm{i})) -U(\tau_\mathrm{i}),  F(\tau_\mathrm{i},
Y_{\mathrm{d}}(\tau_\mathrm{i}))-F(\tau_\mathrm{i},
Y_{\mathrm{u}}(\tau_\mathrm{i})) -U(\tau_\mathrm{i}) \rangle h^2\big]+O(h^{2})
\end{align*}
Hence, $\langle d(\tau_{n^{*}}), d(\tau_{n^{*}}) \rangle \leq
c_{1}h^{2}$ indicating that the controlled trajectory is within the
neighborhood of the desired attractor.

\vspace{10pt}

\textbf{Proof of Theorem \ref{Main_2}}:
Assume that there exists a sufficiently small $h>0$, the interval
$[\tau_{0},\tau^{*}]=\bigcup_{i=0}^{n^{*}}[\tau_\mathrm{i},\tau_\mathrm{i+1}]$,
where $\tau_\mathrm{i+1}=\tau_\mathrm{i}+h$,
$i=1,\cdots,n^{*}$. For
$\tau\in[\tau_\mathrm{i},\tau_\mathrm{i}+h]$, the distance
vector can be defined as
\begin{displaymath}
 d_\mathrm{p}(\tau):=Y_{\mathrm{d}}(\tau_{0})+\int_{\tau_{0}}^{\tau}
 F(\tau_\mathrm{s},Y_{\mathrm{d}}(\tau_\mathrm{s})) \mathrm{d}\tau_\mathrm{s}-
 Y_\mathrm{c}(\tau_{0})-\int_{\tau_{0}}^{\tau}F(\tau_\mathrm{s},Y_\mathrm{u}(\tau_\mathrm{s}),
 u_\mathrm{p}(\tau_\mathrm{s})) \mathrm{d}\tau_\mathrm{s}.
\end{displaymath}
Next, we consider the Taylor expansion of $F(\tau,
Y_\mathrm{u}(\tau),u_\mathrm{p}(\tau))$ within the time interval as
\begin{align*}
&\langle d_\mathrm{p}(\tau_{i+1}), d_\mathrm{p}(\tau_{i+1}) \rangle\\
&\qquad\quad =\;\langle
Y_{\mathrm{d}}(\tau_\mathrm{i})-Y_{\mathrm{u}}(\tau_\mathrm{i})+\int_{\tau_\mathrm{i}}^{\tau_\mathrm{i+1}}
(F(\tau_\mathrm{s},Y_{\mathrm{d}}(\tau_\mathrm{s}))
-F(\tau_\mathrm{i},Y_{\mathrm{u}}(\tau_\mathrm{i}),u_\mathrm{p}(\tau_\mathrm{i}))
-\frac{\mathrm{D} F(\tau_\mathrm{i},Y_{\mathrm{u}}(\tau_\mathrm{i}),u_\mathrm{p}(\tau_\mathrm{i}))}{\mathrm{D} \tau}
(\tau_\mathrm{s}-\tau_\mathrm{i})\nonumber\\
&\qquad\quad\quad-\frac{\mathrm{D} F(\tau_\mathrm{i},Y_{\mathrm{u}}(\tau_\mathrm{i}),
u_\mathrm{p}(\tau_\mathrm{i}))}{\mathrm{D} Y} \dot{Y}_{\mathrm{u}}(\tau_\mathrm{i})
 (\tau_\mathrm{s}-\tau_\mathrm{i})-\frac{\mathrm{D} F(\tau_\mathrm{i},
 Y_{\mathrm{u}}(\tau_\mathrm{i}),u_\mathrm{p}(\tau_\mathrm{i}))}{\mathrm{D} u_\mathrm{p}}
 \dot{u}_\mathrm{p}(\tau_\mathrm{i})(\tau_\mathrm{s}-\tau_\mathrm{i}))\mathrm{d} \tau_\mathrm{s},\nonumber\\
&\qquad\quad\quad
Y_{\mathrm{d}}(\tau_\mathrm{i})-Y_{\mathrm{u}}(\tau_\mathrm{i})+\int_{\tau_\mathrm{i}}^{\tau_\mathrm{i+1}}
(F(\tau_\mathrm{s},Y_{\mathrm{d}}(\tau_\mathrm{s}))
-F(\tau_\mathrm{i},Y_{\mathrm{u}}(\tau_\mathrm{i}),u_\mathrm{p}(\tau_\mathrm{i}))
-\frac{\mathrm{D} F(\tau_\mathrm{i},Y_{\mathrm{u}}(\tau_\mathrm{i}),u_\mathrm{p}(\tau_\mathrm{i}))}{\mathrm{D} \tau}
(\tau_\mathrm{s}-\tau_\mathrm{i})\nonumber\\
&\qquad\quad\quad-\frac{\mathrm{D}
F(\tau_\mathrm{i},Y_{\mathrm{u}}(\tau_\mathrm{i}),u_\mathrm{p}(\tau_\mathrm{i}))}{\mathrm{D}
Y} \dot{Y}_{\mathrm{u}}(\tau_\mathrm{i})
(\tau_\mathrm{s}-\tau_\mathrm{i})-\frac{\mathrm{D}
F(\tau_\mathrm{i},Y_{\mathrm{u}}(\tau_\mathrm{i}),u_\mathrm{p}(\tau_\mathrm{i}))}{\mathrm{D}
u_\mathrm{p}} \dot{u}_\mathrm{p}(\tau_\mathrm{i})
(\tau_\mathrm{s}-\tau_\mathrm{i}))\mathrm{d}\tau_\mathrm{s}\rangle.\end{align*}
By following the same procedure in \textbf{Theorem \ref{Main_1}},
this theorem can be proved.

\bibliographystyle{ieeetr}
\bibliography{bibl}

\end{document}